\documentclass[lettersize,journal]{IEEEtran}
\usepackage{amsmath,amsfonts}
\usepackage{algorithmic}
\usepackage{algorithm}
\usepackage{array}
\usepackage[caption=false,font=normalsize,labelfont=sf,textfont=sf]{subfig}
\usepackage{textcomp}
\usepackage{stfloats}
\usepackage{url}
\usepackage{verbatim}
\usepackage{graphicx}
\usepackage{cite}
\usepackage{amssymb}
\usepackage{subfig}
\usepackage{soul,color}

\begin{document}

\title{A Stable Weighted Residual Finite Element Formulation for the 
Simulation of Linear Moving Conductor Problems}

\author{\IEEEauthorblockN{
    Sethupathy Subramanian and
    Sujata Bhowmick 
    \textit{Member, IEEE}
}
   \thanks{© 2022 IEEE. Personal use of this material is permitted.  Permission from IEEE must be obtained for all other uses, in any current or future media, including reprinting/republishing this material for advertising or promotional purposes, creating new collective works, for resale or redistribution to servers or lists, or reuse of any copyrighted component of this work in other works.}}


\markboth{IEEE Journal on Multiscale and Multiphysics Computational Techniques$~~$-$~~$Accepted article $~~~~~~~~$ DOI 10.1109/JMMCT.2022.3202913}
{Shell \MakeLowercase{\textit{et al.}}: A Sample Article Using IEEEtran.cls for IEEE Journals}

\maketitle

\begin{abstract}
The finite element method is one of the widely employed numerical techniques
in electrical engineering for the study of electric and magnetic fields. 
When applied to the moving conductor problems, the finite element method
is known to have numerical oscillations in the solution. To resolve this,
the upwinding techniques, which are developed for the transport equation 
are borrowed and directly employed for the magnetic induction equation. 
In this work, an alternative weighted residual formulation is explored 
for the simulation of the linear moving conductor problems. 
The formulation is parameter-free and the stability of the formulation 
is analytically studied for the 1D version of the moving conductor problem.
Then the rate of convergence and the accuracy are illustrated with the help of
several test cases in 1D as well as 2D. Subsequently, the stability of the 
formulation is demonstrated with a 3D moving conductor simulation.

\end{abstract}


\begin{IEEEkeywords}
Numerical Stability, Parameter-free, Moving Conductor, Advection-Diffusion
\end{IEEEkeywords}



\section{Introduction} \label{sec:1}


The finite element method is a widely used numerical technique for the
design and analysis of electrical machines and instruments. The Galerkin
finite element (GFE) formulation is known to produce highly accurate
solutions for electrostatic and magnetostatic simulations. However, 
for the simulation of linear moving conductors, such as linear induction motor, 
magnetic brakes, electromagnetic flowmeter etc., the numerical stability
problem occurs at high velocities \cite{mc3eb1, mc1lim1, su1}.

The numerical stability problem at high velocities is also common to 
the transport equation of fluid dynamics \cite{cdbook}. 
The numerical oscillation in the simulation of transport equation is 
well studied over several decades \cite{quada1, revfic1}.
It is observed from the finite difference formulation that the stability is
restored, when the central difference of the first order term is
reinforced/replaced with the one-sided difference along the flow
direction \cite{upfdm1_sp, upfdm1_ru}. Following this observation, 
the common technique in finite element method is to upwind the weight
function along the flow direction,
so to make the first order term more one-sided \cite{up1, supg1}.

These upwinding techniques are borrowed from the transport equation and 
directly employed for the moving conductor problems 
\cite{mc6tf1, mc2av1, mc4ge1, mc5mc1, mcsupg_mcfit, mcsupg_cable, mcsupg_mfluid}. 
Even though the upwinding techniques solve
the numerical instability problem, there are other issues associated with the
application of upwinding schemes, such as, crosswind diffusion and
 erroneous solution at the transverse boundary \cite{supg1, discop1, discop}. 
 These errors were first observed for the transport equation and several 
 remedies have been suggested, with partial success 
 \cite{soldreview1, soldreview2}. The error at the transverse boundary
 observed for the moving conductor problems as well; and a solution is
 suggested in a recent literature \cite{nemosu, sus3}.
 
Given these, in this work, a weighted residual finite element formulation 
is suggested, by eliminating the first order derivatives in the 
governing equation of the moving conductor problem, in a way it is consistent. 
Thus, the term which introduces numerical instability is eliminated and
numerical stability can be achieved. It can also be noted that such a formulation
does not require a \emph{stabilization parameter - $\tau$}, like the 
upwinding schemes. Thus, this formulation remains parameter-free. 

Firstly, stability analysis is performed for the weighted residual 
formulation in 1D. Then the numerical exercises are carried out for 1D, as well 
as, 2D problems to check the accuracy and the convergence. In addition to this, 
a test case involving the magnetic material $\mu_r > 1$ is simulated in 2D 
as well as in 3D, to verify the accuracy of the formulation with multiple
materials.  In the next section the analysis for the 1D problem is discussed.


\section{Analysis on the 1D moving conductor problem} \label{sec:2}


The 1D version of the moving conductor problem is derived in \cite{su1} and
it can be described as follows. In this, the source magnetic
field is applied perpendicular to the plane of the paper ($x$-direction).
The conductor of an infinite dimension,
is moving along the horizontal $z$-axis with velocity $u_z$,
permiability $\mu$ and conductivity $\sigma$.
The reaction magnetic field $b_x$ is arising out of the motion
and its vector potential is $A_y$. The governing equation for
this problem can be written as,
\begin{align} \label{eq:1Dge}
    -\dfrac{d^2 A_y}{dz^2} + \mu \sigma u_z \dfrac{d A_y}{dz} = 
                             \mu \sigma u_z B_x
\end{align}
In this, the first derivative term can be eliminated from the governing equation, by
writing the advective-first order term as $b_x = -{dA_y}/{dz}$. The 
equation (\ref{eq:1Dge}) becomes,
\begin{align} \label{eq:1Dmc_1a}
    -\dfrac{d^2 A_y}{dz^2} - \mu \sigma u_z b_x = 
                             \mu \sigma u_z B_x
\end{align}
The Galerkin finite element formulation of (\ref{eq:1Dmc_1a}) is
written below with $N$ as weight (shape) function and the integration-by-parts is
applied to the second-derivative diffusion term. In this, $A_y$, $b_x$ and the 
weight function $N$ belong to $H^1$ function space.

\begin{align} \label{eq:wr1Dmc_1}
    \int_\Omega \dfrac{dN}{dz}\dfrac{d A_y}{dz} ~d\Omega - 
                          \mu \sigma u_z \int_\Omega N b_x ~d\Omega = 
                          \mu \sigma u_z \int_\Omega N B_x ~d\Omega 
\end{align}
Now the first derivative term is replaced with $b_x$ and it becomes a new
unknown along with $A_y$. Therefore, we need a second equation to solve
this system. The second equation is written below with the weight function of 
$dN/dz$, so as to eliminate the numerically unstable $N (dA_y/dz)$ term.

\begin{align} \label{eq:wr1Dmc_2}
    \int_\Omega \dfrac{dN}{dz} \dfrac{d b_x}{dz} ~d\Omega + 
     \mu \sigma u_z \int_\Omega \dfrac{dN}{dz} \dfrac{d A_y}{dz} ~d\Omega 
     ~\dots \nonumber \\ 
     = \mu \sigma u_z \int_\Omega \dfrac{dN}{dz} B_x ~d\Omega 
\end{align}
In this, the second derivative term $-{d^2 A_y}/{dz^2}$ of (\ref{eq:1Dge}) 
is replaced with ${d b_x}/{dz}$. 
Equations (\ref{eq:wr1Dmc_1}) and (\ref{eq:wr1Dmc_2}) together can form a 
stable weighted formulation without any upwinding. In the next subsection,
the stability analysis is carried out for (\ref{eq:wr1Dmc_1}) and 
(\ref{eq:wr1Dmc_2}).


\subsection{Stability analysis for the 1D weighted residual formulation}


It can be readily noted that the governing equation of the 1D moving conductor
problem (\ref{eq:1Dge}) is same as that of the transport equation. 
However, in moving conductor problems, the boundary conditions on the
magnetic vector potential $\bf{A}$ are not hard-imposed. 
The boundary conditions are always `${\bf A} = 0$' very far from the source 
magnetic field $\bf{B^a}$ or natural boundary conditions are chosen, if the 
problem permits \cite{su1}. Under this circumstance,
the $Z$-transform based $pole$, $zero$ analysis is a good indicator for
the numerical stability \cite{su1, su2} in moving conductor problems.
Therefore, the same is employed here. 

The difference equation form of the finite element equation (\ref{eq:wr1Dmc_1}),
with linear elements, at the $n^{th}$ node can be written as,
\begin{align} \label{eq:wr1Dmc_d1}
    - A_{y[n-1]} + 2 A_{y[n]} - A_{y[n+1]}   
    - \dfrac{Pe}{3} ( b_{x[n-1]}
    ~\dots \nonumber \\
  + 4 b_{x[n]} + b_{x[n+1]} ) 
    = \dfrac{Pe}{3} ( B_{x[n-1]} + 4 B_{x[n]} + B_{x[n+1]} )
\end{align}
and similarly the difference equation form of the finite element equation 
(\ref{eq:wr1Dmc_2}), with linear elements, at the $n^{th}$ node can be written as,
\begin{align} \label{eq:wr1Dmc_d2}
          - b_{x[n-1]} + 2 b_{x[n]} - b_{x[n+1]}
 + 2 Pe ( - A_{y[n-1]}
    ~\dots \nonumber \\
    + 2 A_{y[n]} - A_{y[n+1]} ) = Pe (B_{x[n-1]} - B_{x[n+1]})
\end{align}
where $Pe$ is the Peclet number and it is defined as 
$Pe = \mu \sigma u_z \Delta z / 2$. The Peclet number indicates the 
relative strength of advection over the diffusion in the difference 
equation.
Taking $Z$-transform of the difference equations (\ref{eq:wr1Dmc_d1}) and
(\ref{eq:wr1Dmc_d2}),
\begin{align} \label{eq:wr1Dmc_z1}
    - A_y ( Z^2 - 2 Z + 1) - \dfrac{Pe}{3} b_x ( Z^2 + 4 Z + 1 ) 
    ~\dots \nonumber \\
                           = \dfrac{Pe}{3} B_x ( Z^2 + 4 Z + 1 )
\end{align}
\begin{align} \label{eq:wr1Dmc_z2}
    b_x ( Z^2 - 2 Z + 1 ) + 2 Pe A_y ( Z^2 - 2 Z + 1 )
                                = Pe B_x (Z^2 - 1)
\end{align}
Now, substituting the expression for $b_x$ from (\ref{eq:wr1Dmc_z1}), in
(\ref{eq:wr1Dmc_z2}) and then approximating for $Pe >> 1$, the final 
transfer function between $A_y$ and $B_x$ can be written as,
\begin{align} \label{eq:wr1Dmc_tf}
    \dfrac{A_y}{B_x} \approx \dfrac{Z^2 - 1}{2 (Z^2 - 2 Z + 1)}
\end{align}
The \emph{poles} of the transfer function (\ref{eq:wr1Dmc_tf}) are $Z = 1, 1$
and they are positive; indicating a stable formulation. In other words,
the roots of the difference equations are positive for $Pe >> 1$, resulting
in a non-oscillatory system for the moving conductor problems. In the next
subsection, the simulation results from 1D moving conductor problem
are presented.


\subsection{Simulation results from 1D moving conductor problem}

\begin{figure}
		\centering
		\mbox{\subfloat[]{\label{f:1Dmc_Pe4}
		\includegraphics[scale=0.45]{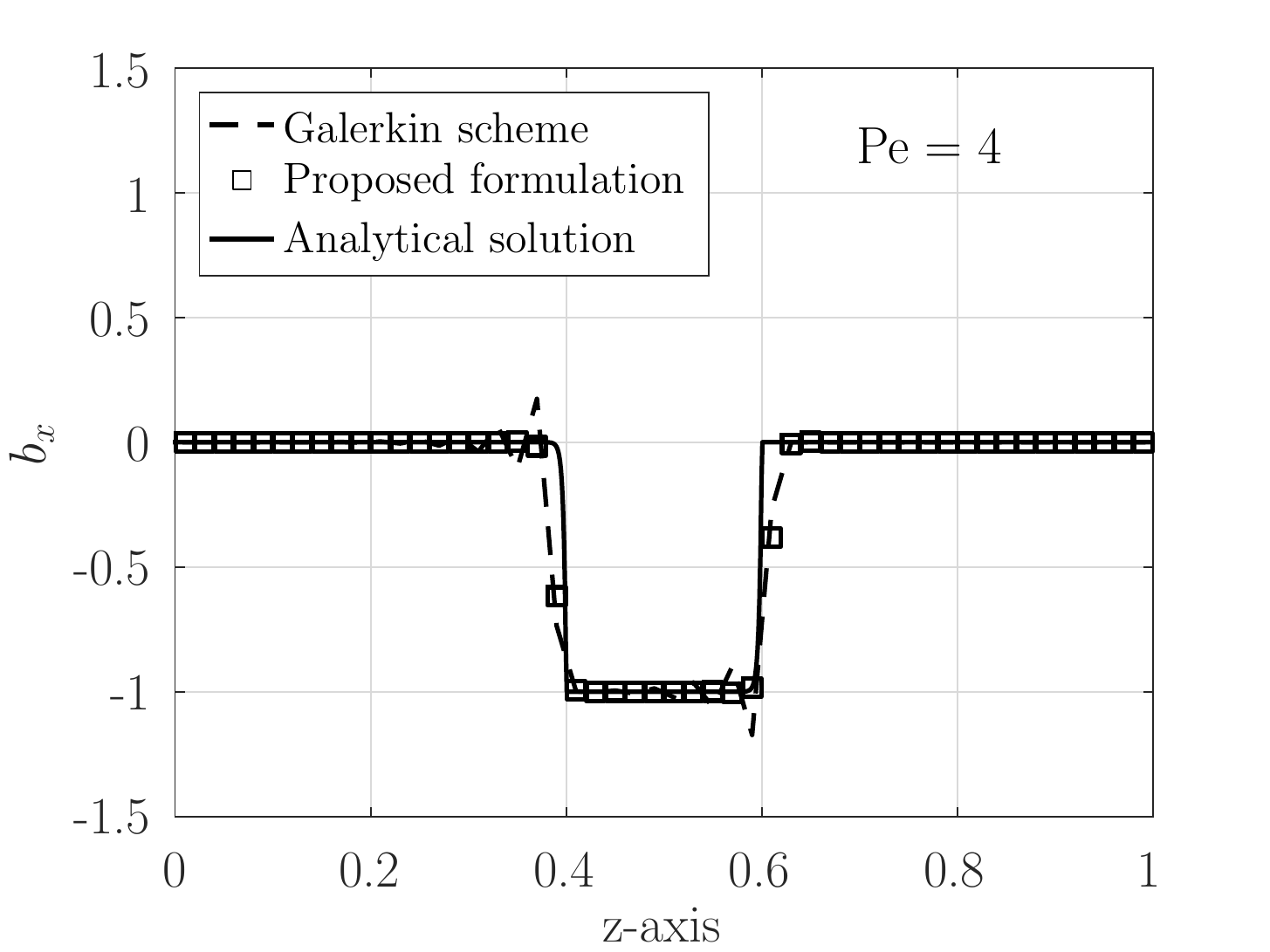}}}
		\mbox{\subfloat[]{\label{f:1Dmc_Pe400}
		\includegraphics[scale=0.45]{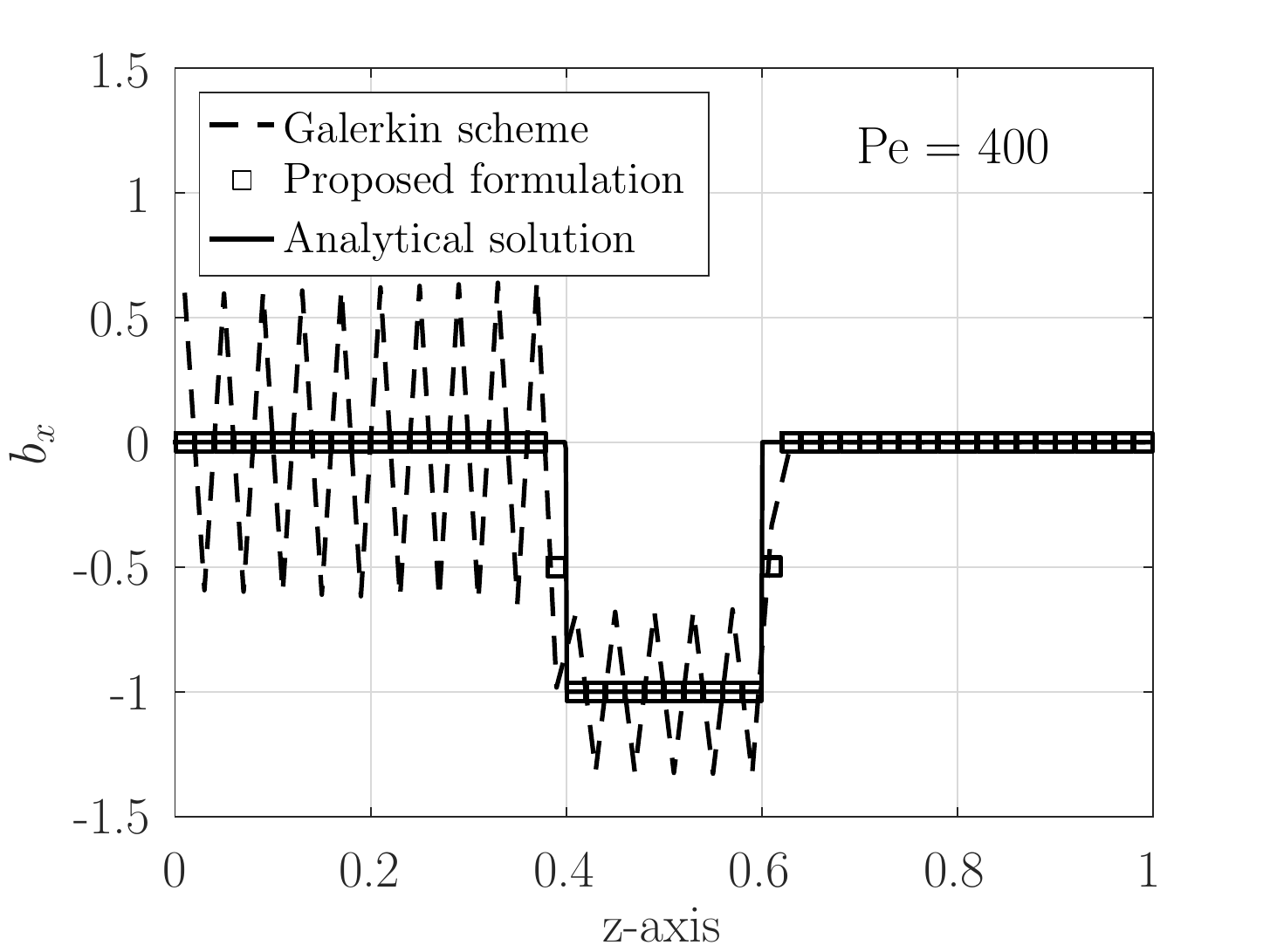}}}

		\caption{ Comparison of simulation results for the 1D moving
        conductor problem (a) Pe = 4, (b) Pe = 400.}
		\label{f:1Dmc}
\end{figure}

The finite element simulation of (\ref{eq:1Dge}) is carried out, 
using the Galerkin formulation and the weighted residual formulation, which
is described by equations (\ref{eq:wr1Dmc_1}) and (\ref{eq:wr1Dmc_2}). 
The simulation domain spans $0 \leq z \leq 1$ and the input magnetic field
$B_x = B, for ~ 0.4 \leq z \leq 0.6$. The simulations are carried out for
different set of parameters and the results are observed to be accurate and 
stable.  The sample simulation results for $Pe=4$ and $Pe=400$ are 
presented in Fig. \ref{f:1Dmc}. The figure \ref{f:1Dmc} shows the reaction
magnetic field $b_x = -dA_y/dz$ obtained from the i) proposed weighted residual
formulation, ii) Galerkin scheme and iii) analytical solution described 
in \cite{su1}.
\begin{table}[h!]
\begin{center}
\caption{ Measured values of error in $b_x = -dA_y/dz$ (first derivative)\
          with the proposed formulation for 1D moving conductor problem
          with $\mu \sigma u_z = 1000$ }
\label{t:1Dmc}
\begin{tabular}{|c|c|c|c|}
\hline
 Number of  &    L2     & Absolute  & Expt. Order of       \\
 Elements   &    error  & error     & Convergence          \\
\hline
  50        & 2.49e-03  & 2.04e-02  &               \\
  100       & 1.76e-03  & 1.05e-02  & 0.96          \\
  200       & 1.21e-03  & 5.40e-03  & 0.95          \\
  400       & 7.84e-04  & 2.64e-03  & 1.03          \\
  800       & 4.50e-04  & 1.29e-03  & 1.04          \\
\hline
\end{tabular}
\end{center}
\end{table}

In table \ref{t:1Dmc}, the absolute and the {root mean-squared (rms - L2)} 
values of errors are displayed along with the experimental order of convergence. 
The errors are calculated for the reaction magnetic field $b_x = -dA_y/dz$, 
which is a quantity of interest in the moving conductor simulations.
It can be seen that the weighted residual formulation provides stable as well 
as converging results for the 1D case.

The stability of the weighted residual formulation of 
(\ref{eq:wr1Dmc_1}), (\ref{eq:wr1Dmc_2}) are due to the non-oscillatory
\emph{poles/roots} present in the difference equation. This is in 
contrast to the stable formulations presented in \cite{su1, su2},
where the stability of the formulation is achieved by canceling the 
oscillatory \emph{poles} with the help of \emph{zeros} of the input 
magnetic field. Therefore, the formulations of \cite{su1, su2} depends 
on the representation of the input magnetic field for the stability,
while the formulation presented in (\ref{eq:wr1Dmc_1}), (\ref{eq:wr1Dmc_2})
are not dependent on the input magnetic field. Given this situation, it
may be worthwhile to check the stability of the formulation for the
transport equation. This is carried out in the next subsection.


\subsection{Simulation results from 1D transport equation}


The transport equation describes the transport of a physical variable $\psi$
by means of advection and diffusion. The amount of diffusion is described
by the diffusivity parameter $k$ and $u$ is the velocity of motion.
The transport equation with the source term $S$ is given by,
\begin{align} \label{eq:1Dtr}
    -k \dfrac{d^2 \psi}{dz^2} + u \dfrac{d \psi}{dz} = S
\end{align}
The weighted residual formulation can be written as,
\begin{align} \label{eq:wr1Dtr_1}
    \int_\Omega \dfrac{dN}{dz}\dfrac{d \psi}{dz} ~d\Omega + 
                          \dfrac{u}{k} \int_\Omega N F_z ~d\Omega = 
                          \dfrac{1}{k} \int_\Omega N S   ~d\Omega
\end{align}
\begin{align} \label{eq:wr1Dtr_2}
    -\int_\Omega \dfrac{dN}{dz} \dfrac{d F_z}{dz} ~d\Omega +
                \dfrac{u}{k} \int_\Omega \dfrac{dN}{dz} \dfrac{d \psi}{dz} ~d\Omega =
                \dfrac{1}{k} \int_\Omega \dfrac{dN}{dz} S  ~d\Omega
\end{align}
where, $N$ is the shape function and $F_z$ is the flux, defined as
$F_z = d\psi/dz$. For the simulation, two standard test cases are considered.
The first test problem (TP1) has source term $S=0$ and $\psi=0$ at $z=0$ and 
$\psi = 1$ at $z=1$. The $u/k$ ratio is taken to be 400 for this test case.
The second test problem (TP2) has source term $S=z^2$ and $\psi=0$ at $z=0$ 
and $z=1$. The $u/k$ ratio is taken to be 200 for this second test case.

\begin{figure}
		\centering
		\mbox{\subfloat[]{\label{f:1Dsm}
		\includegraphics[scale=0.45]{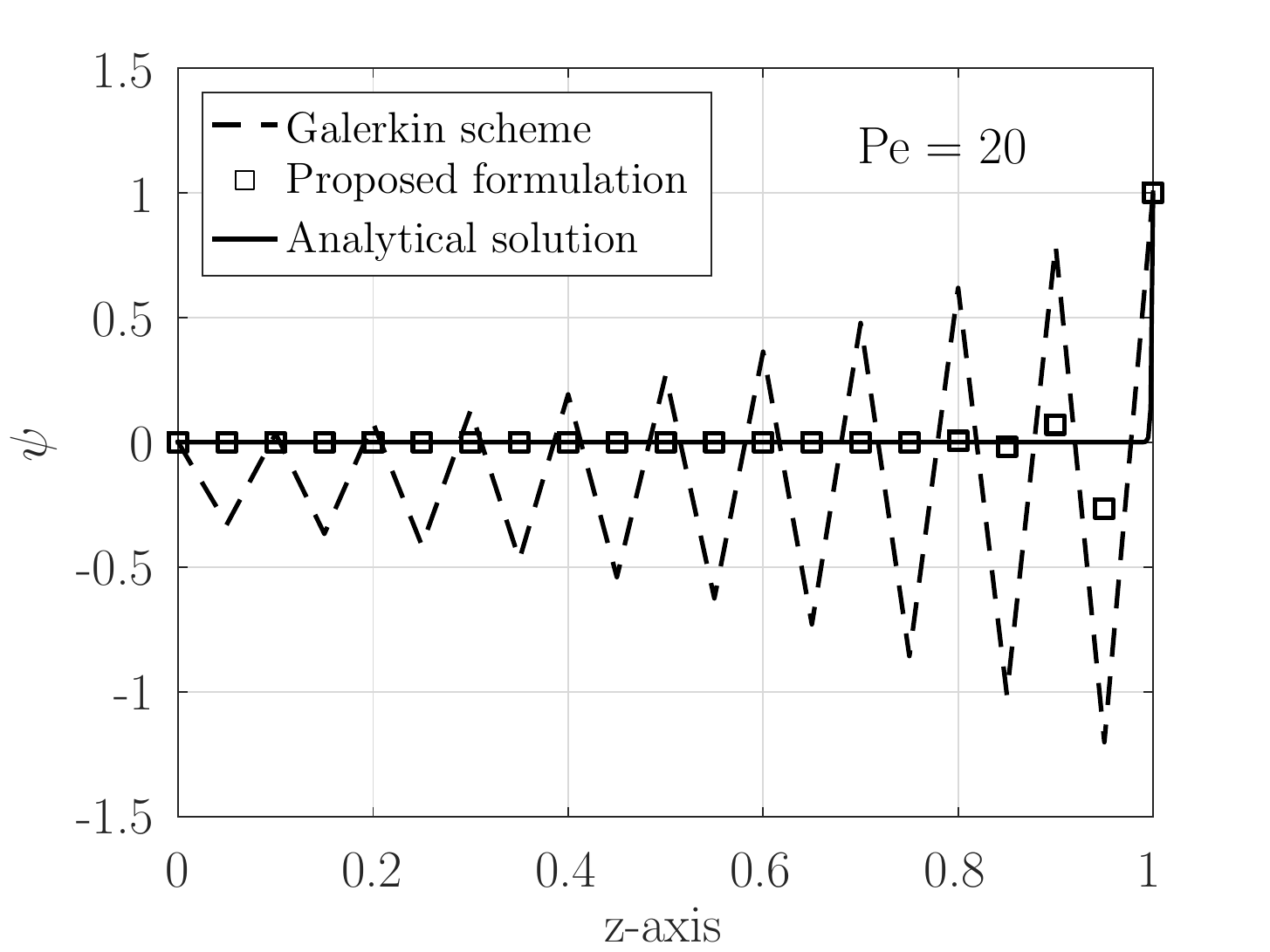}}}
		\mbox{\subfloat[]{\label{f:1Dx2}
		\includegraphics[scale=0.45]{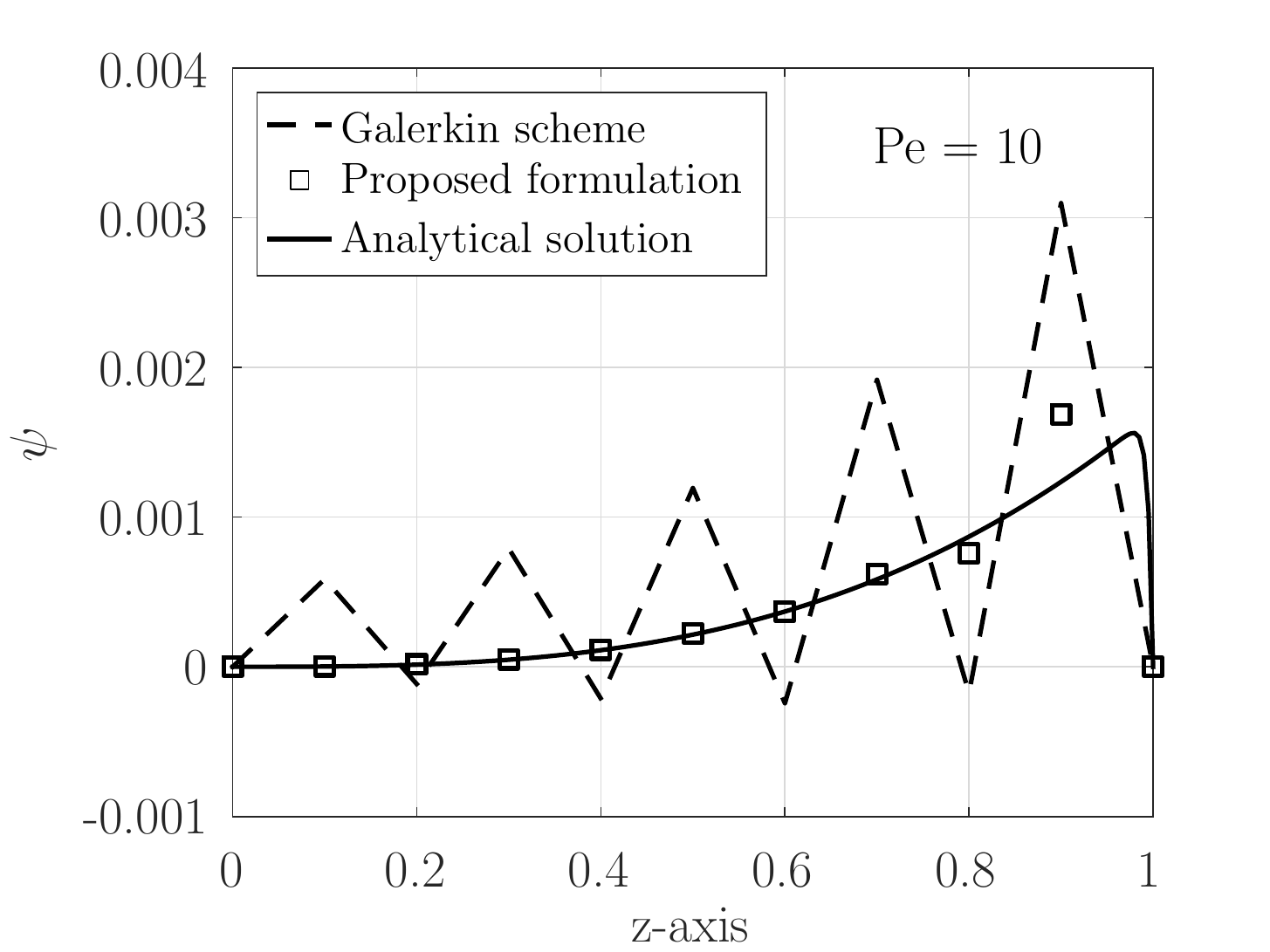}}}
		\caption{ Comparison of simulation results for the 1D
        transport equation 
        (a) TP1 
        (b) TP2 
        }
		\label{f:1Dtr}
\end{figure}

\begin{table}[h!]
\begin{center}
\caption{Error in $\psi$ with the proposed formulation for TP1 and TP2}
\label{t:1Dsmx2}
\begin{tabular}{|c|c|c|c|c|}
\hline
    & Number of  &    L2     & Absolute  & Expt. Order of       \\
    & Elements   &    error  & error     & Convergence          \\
\hline
    &      20   & 4.15e-03 & 2.62e-02    &            \\
    &      40   & 2.53e-03 & 1.16e-02    & 1.18       \\
TP1 &      80   & 1.31e-03 & 4.24e-03    & 1.45       \\
    &     160   & 5.08e-04 & 1.07e-03    & 1.99       \\
    &     320   & 1.55e-04 & 2.44e-04    & 2.13       \\
\hline
    &      10   & 1.40e-05 & 9.16e-05    &            \\
    &      20   & 8.55e-06 & 4.04e-05    & 1.18       \\
TP2 &      40   & 4.42e-06 & 1.47e-05    & 1.46       \\
    &      80   & 1.72e-06 & 3.74e-06    & 1.98       \\
    &     160   & 5.23e-07 & 8.59e-07    & 2.12       \\
\hline
\end{tabular}
\end{center}
\end{table}

The simulation results are displayed in figure \ref{f:1Dtr}. The results show 
that the proposed formulation performs better than the Galerkin scheme. 
However, the absolute stability observed in the moving conductor problem
(see figure \ref{f:1Dmc}) is not observed for the transport equation.
This can be attributed to the hard boundary conditions set at either
end of the simulation, which introduces a steep slope near the boundary.

The table \ref{t:1Dsmx2}, displays the absolute and rms errors
obtained from the weighted residual formulation for the first and second cases. 
The errors are calculated for the variable $\psi$, by 
comparing it with the respective analytical solution. The results show the 
expected experimental order of convergence for the 
weighted residual formulation.

Thus, even though the proposed formulation did not perform to 100\%,
the overall accuracy is achieved for the transport equation. 
{It may be noted that for the 1D problems discussed so far, the upwinding
formulations, such as the Streamline upwinding/Petrov-Galerkin (SU/PG)
scheme would provide a highly accurate solution; close to the
analytical solution. This is due to the fact that the upwinding stabilization
parameter is derived by matching the numerical formulation with the
analytical solution} \cite{quada1}.
In the next section, analysis for the 2D moving conductor problems are
described.


\section{Analysis with 2D moving conductor problems} \label{sec:3}


The moving conductor problems, when described in two dimensions can exhibit
two kinds of circulations; i) circulation of vector potential or current
ii) circulation of magnetic field.
This is due to the curl nature of the governing equation and it can be
written as follows \cite{su1, sus3}:
\begin{equation} \label{eqge1_r}
\sigma \nabla \phi ~-~  (\nabla \cdot \dfrac{1}{\mu} \nabla) {\bf{A}} - \sigma~ {\bf{u}} \times \nabla \times {\bf{A}} = \sigma~ {\bf{u}} \times {\bf{B_{a}}}
\end{equation}
\begin{equation} \label{eqge2_r}
\begin{split}
\nabla \cdot (\sigma \nabla \phi) - \nabla \cdot (\sigma ~ {\bf{u}} \times \nabla \times {\bf{A}}) =  \nabla \cdot (\sigma ~ \bf{u} \times \bf{B_{a}})
\end{split}
\end{equation}
where, ${\bf u}$ is the velocity and $\phi$ is the electric scalar potential.
In next subsection, the weighted residual formulation for a problem containing 
the circulation of vector potential or current is described.


\subsection{2D problem with the circulation of {\bf A}}


\begin{figure}
        \centering
        \includegraphics[scale=0.3]{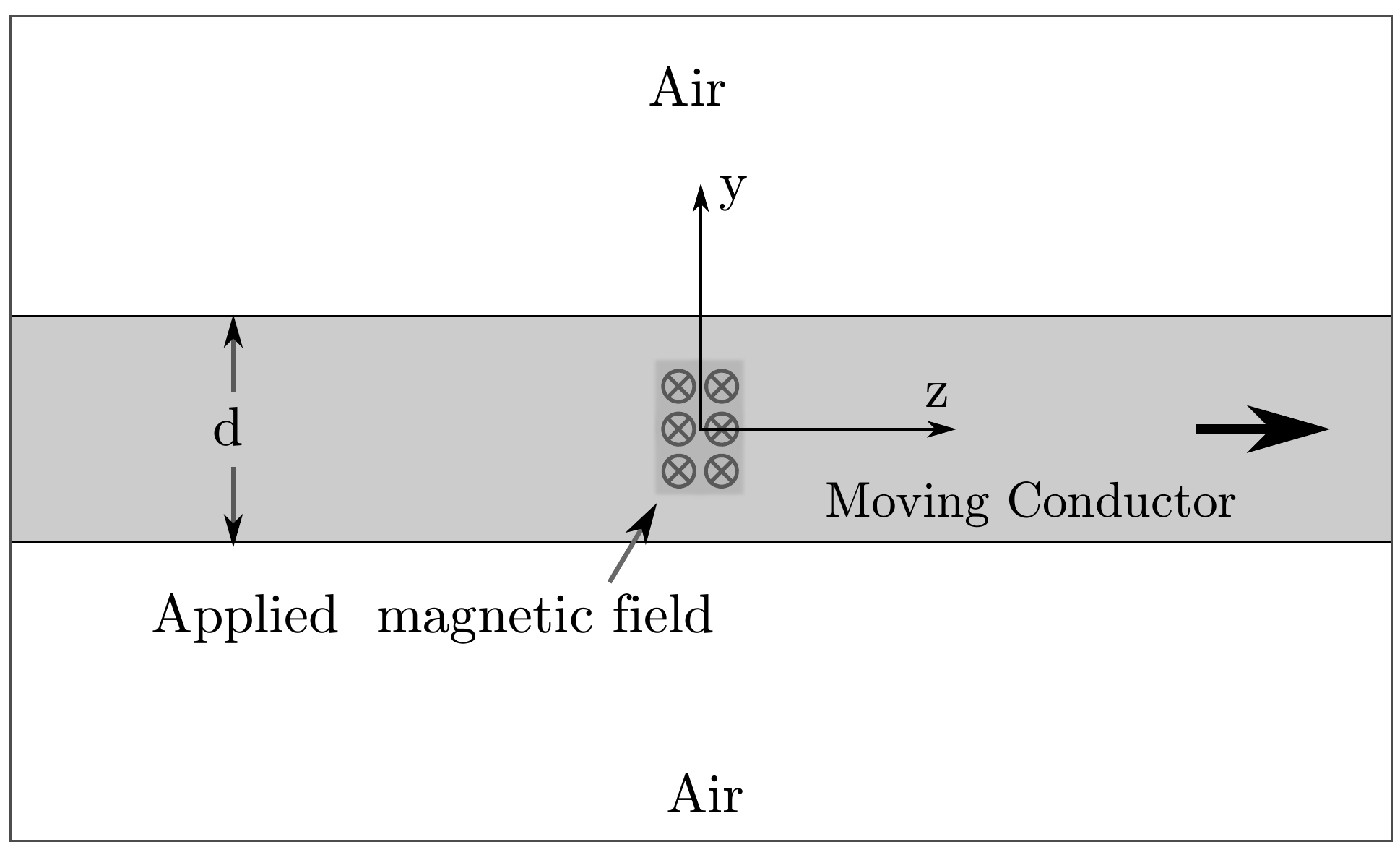}
        \caption{ Schematic of the 2D moving conductor problem with
        the circulation of $\bf A$}
        \label{f:2Dmc_A}
\end{figure}

The schematic of the 2D problem is shown in Fig. \ref{f:2Dmc_A}. The problem is
defined along the $zy-$axis. The input magnetic field is directed perpendicular 
to the plane and it is directed along the x-axis ($B_x^a$). The conductor
of width $d$ is moving along the $z$-axis with the velocity $u_z$. The simulation
boundaries are chosen to be far from the input magnetic field.

In this problem, conductor motion creates a reaction magnetic field $b_x$,
which tries to cancel the input magnetic field following the Lenz's law.
Therefore the current circulation and the circulation of the magnetic vector
potential ${\bf A}$ are along the $zy-$plane and it has two components,
$A_y$ and $A_z$. 
The governing equations of this problem can be written as \cite{su2},

\begin{equation} \label{eq:ge2dphi_A}
    \sigma \nabla^2 \phi + 
                    \sigma u_z \dfrac{\partial^{2}A_y}{\partial y \partial z} - 
                    \sigma u_z \dfrac{\partial^{2}A_z}{\partial y^2         } = 
                    \sigma u_z \dfrac{\partial B_x   }{\partial y}
\end{equation}

\begin{equation} \label{eq:ge2day_A}
 \sigma \dfrac{\partial \phi}{\partial y} - \dfrac{1}{\mu}\nabla^2 A_y 
    + \sigma u_z \dfrac{dA_y}{dz} 
    - \sigma u_z \dfrac{dA_z}{dy} = \sigma u_z B_x
\end{equation}

\begin{equation} \label{eq:ge2daz_A}
 \sigma \dfrac{\partial \phi}{\partial z} - \dfrac{1}{\mu}\nabla^2 A_z = 0
\end{equation}

Here, the first derivative term is arising from 
${\bf u} \times \nabla \times {\bf A}$. This term can be rewritten as 
${\bf u} \times {\bf b}$, following the 1D problem defined in equations
(\ref{eq:wr1Dmc_1}), (\ref{eq:wr1Dmc_2}). In this problem, $\bf A$ has the 
circulation and only the x-component of the reaction magnetic field $b_x$
exists. In other words, the equations (\ref{eq:ge2day_A}),
(\ref{eq:ge2daz_A}) when written in terms of $b_x$, they become 2 equations and 
one unknown - $b_x$. Therefore, this is taken into account in the formulation
of the weighted residual scheme.

The Galerkin finite element formulation of the equations (\ref{eq:ge2dphi_A}), 
(\ref{eq:ge2day_A}), (\ref{eq:ge2daz_A}) are
written below with $N$ as weight (shape) function and the integration-by-parts is
applied to second-derivative terms. Also, similar to the 1D case, the 
${\bf u} \times \nabla \times {\bf A}$ is replaced with the ${\bf u} \times {\bf b}$
term.

\begin{equation} \label{eq:wr2dphi_A}
    \int_\Omega \sigma \nabla N \cdot \nabla \phi ~ d\Omega -
    \int_\Omega \sigma u_z \dfrac{\partial N}{\partial y} b_x ~d\Omega=
    \int_\Omega \sigma u_z \dfrac{\partial N}{\partial y} B_x^a ~d\Omega
\end{equation}

\begin{align} \label{eq:wr2day_A}
\int_\Omega \sigma N \dfrac{\partial \phi}{\partial y} ~d\Omega 
    + \dfrac{1}{\mu} \int_\Omega \nabla N \cdot \nabla A_y ~d\Omega 
    ~\dots \nonumber \\
    - \int_\Omega \sigma u_z N b_x ~d\Omega 
    = \int_\Omega \sigma u_z N B_x^a ~d\Omega
\end{align}

\begin{equation} \label{eq:wr2daz_A}
\int_\Omega \sigma N \dfrac{\partial \phi}{\partial z} ~d\Omega 
    + \dfrac{1}{\mu} \int_\Omega \nabla N \cdot \nabla A_z ~d\Omega = 0
\end{equation}

Now, following the 1D formulation (\ref{eq:wr1Dmc_2}), one can write two more 
equations from (\ref{eq:ge2day_A}), (\ref{eq:ge2daz_A}) with the weight functions
$\partial N/\partial z$, $-\partial N/\partial y$ respectively. The
weight functions are chosen such that the formulation has \emph{bilinear form}
in either $\bf A$ or $\bf b$ or both. 
%
%
%
%
The two equations are then summed and form a single equation for the one
extra unknown - $b_x$; it is written as,
\begin{align} \label{eq:wr2day_A_2}
    \int_\Omega \sigma \dfrac{\partial N   }{\partial z} 
                       \dfrac{\partial \phi}{\partial y} ~d\Omega 
  - \int_\Omega \sigma \dfrac{\partial N   }{\partial y} 
                       \dfrac{\partial \phi}{\partial z} ~d\Omega 
        ~ \dots \nonumber \\
    + \dfrac{1}{\mu} \int_\Omega \nabla N \cdot \nabla b_x ~d\Omega
        - \int_\Omega \sigma u_z \dfrac{\partial N   }{\partial z}
                                 \nabla \times {\bf A} \cdot \hat{x} ~d\Omega
        ~ \dots \nonumber \\
        = \int_\Omega \sigma u_z \dfrac{\partial N   }{\partial z} B_x^a ~d\Omega
\end{align}
\begin{figure*}
		\centering
		\mbox{\subfloat[]{\label{f:2DmcA_GL}
		\includegraphics[scale=0.42]{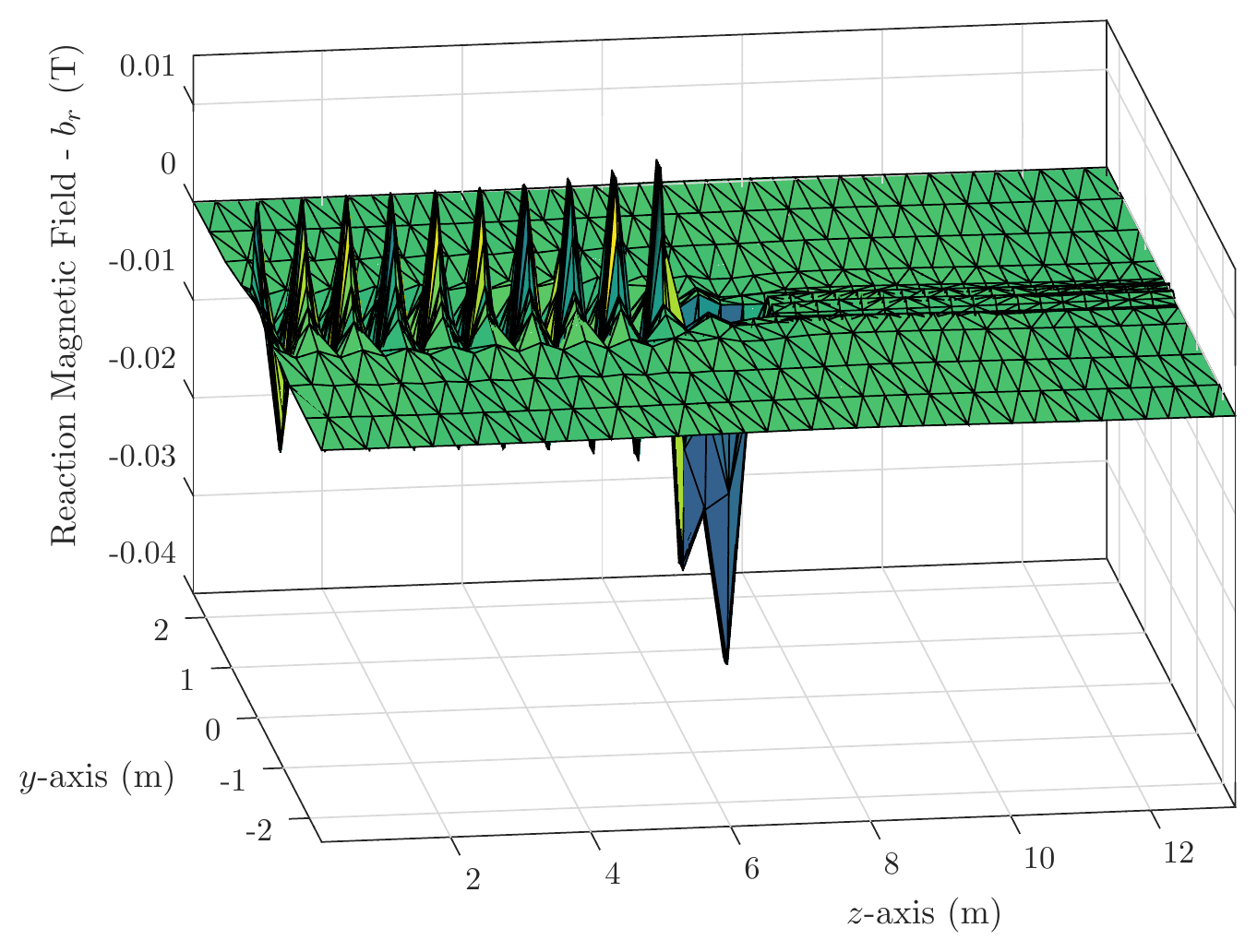}}}
		\mbox{\subfloat[]{\label{f:2DmcA_SUPG}
		\includegraphics[scale=0.42]{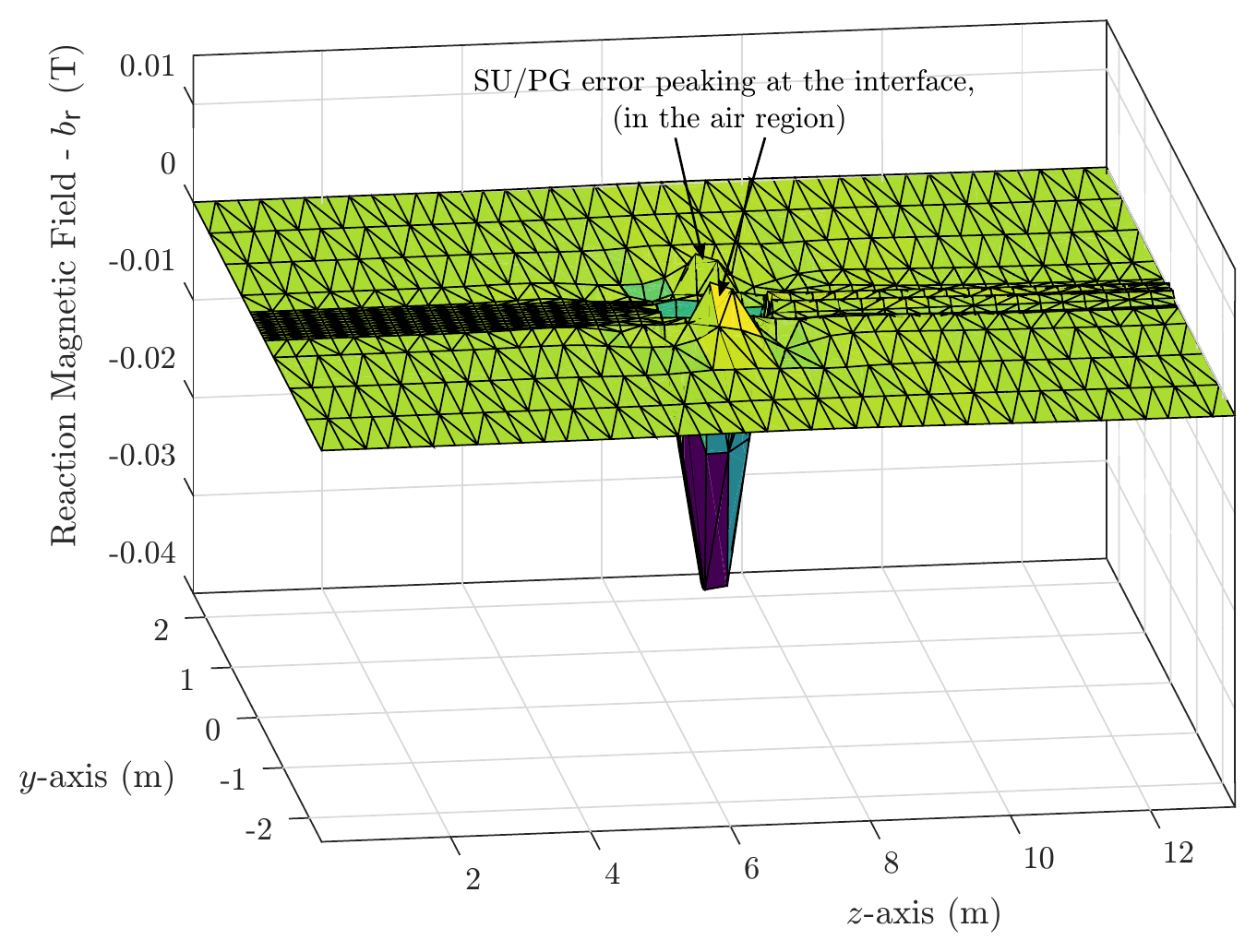}}}
		\mbox{\subfloat[]{\label{f:2DmcA_WR}
		\includegraphics[scale=0.42]{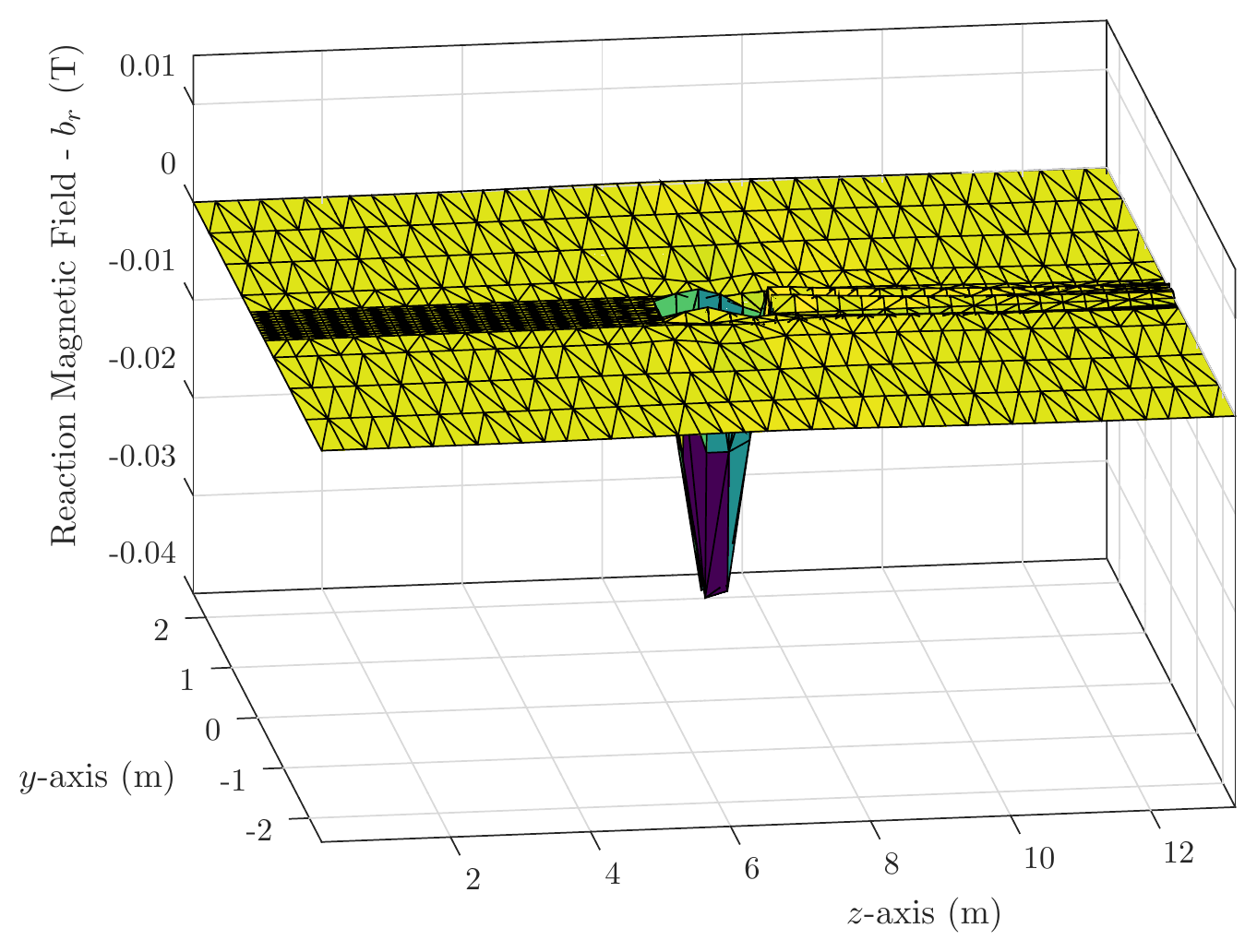}}}
		\caption{ Reaction magentic field $b_x = dA_z/dy - dA_y/dz$ 
        for the 2D moving conductor problem with the circulation of $\bf A$
        (a) Galerkin scheme (b) SU/PG scheme \cite{supg1, nemosu, sus3}
        (c) Proposed formulation}
		\label{f:2DmcA}
\end{figure*}
\begin{table}[h!]
\begin{center}
\caption{ Measured values of error in the reaction magnetic field
          $b_x = dA_z/dy - dA_y/dz$ with the proposed formulation
          for the 2D moving conductor problem with
          the circulation of $\bf A$}
\label{t:2DmcA}
\begin{tabular}{|c|c|c|c|}
\hline
 Number of & L2       & Absolute    & Expt. Order of \\
 Elements  & error    & error       & Convergence\\
\hline
   10240   & 4.52e-06 & 4.44e-04    &            \\
   20480   & 2.85e-06 & 2.02e-04    & 1.14       \\
   40960   & 1.81e-06 & 9.83e-05    & 1.04       \\
   81920   & 1.01e-06 & 4.41e-05    & 1.16       \\
\hline
\end{tabular}
\end{center}
\end{table}
The equations (\ref{eq:wr2dphi_A}), (\ref{eq:wr2day_A}), (\ref{eq:wr2daz_A}), 
(\ref{eq:wr2day_A_2}) form the weighted residual formulation for this 2D case.
Simulations are performed with bilinear quadrilateral elements, 
for different $Pe$ values and stable solutions
are observed. A sample simulation plot is presented for $Pe=1500$ in 
Fig. \ref{f:2DmcA}, where Fig, \ref{f:2DmcA_GL} presents the solution from the
Galerkin scheme and Fig. \ref{f:2DmcA_WR} presents the solution from the
weighted residual formulation. In addition to this, an accuracy study is 
conducted and the results are presented in table \ref{t:2DmcA}. 
In table \ref{t:2DmcA}, the error in $\nabla \times {\bf A}$ is computed
by comparing the solution with the solution obtained from a very
fine discretisation.
The error and the rate of convergence shows that the weighted residual
formulation produces accurate and converging solutions for this
2D case.
{
For quick comparison, Fig. {\ref{f:2DmcA_SUPG}} displays the result obtained 
from the SU/PG scheme. The magnetic field has the error peaking at the material
interface and the peak occurs in the air-region adjacent to the moving conductor.
This would lead to non-physical current circulation in the air-region }
\cite{nemosu, sus3}.

In the next subsection, the second case of the 2D moving conductor is
described; where the circulation of the reaction magnetic field ${\bf b}$ is
present, instead of the circulation of ${\bf A}$. For this, the `Testing
Electromagnetic Analysis Methods' (TEAM) problem No. 9a is chosen \cite{team9},
and it is detailed in the next subsection.


\subsection{2D problem with the circulation of {\bf b}}

A schematic representation of the TEAM-9a problem is provided in
Fig. \ref{team9a_sch}.
The TEAM-9a problem involves an infinite ferromagnetic material with 
the conductivity $\sigma=5\times 10^6~Sm^{-1}$ and the relative permeability of 
cases with $\mu_r=1,50$. The ferromagnetic material has a cylindrical bore of 
radius $r_i = 14mm$. Inside the bore, a concentric current loop of 
diameter $r_c = 12mm$ is carrying a current of $1A$ and it moves at an 
uniform velocity in the bore. 
This is an axisymmetric problem along the $z,~r$ axes, and it has no
variation along the $\theta$-axis.
For the analysis, the highest case with the velocity of $v=100ms^{-1}$ is 
considered. The finite element mesh is denser close to the current loop and 
becomes coarser as moving away from the current loop. The resulting $Pe$ varies 
from $5$ to $200$, due to the varying discretisation.

The coupled governing equation for this axisymmetric problem can be written in 
terms of the magnetic vector potential $A_\theta$
and the radial reaction magnetic field $b_r$ as,

\begin{equation}\label{eq:2Dt9ge}
 -\dfrac{1}{\mu}\left( \dfrac{\partial^2 A_\theta}{\partial r^2} + 
                       \dfrac{\partial^2 A_\theta}{\partial z^2} +
           \dfrac{1}{r}\dfrac{\partial   A_\theta}{\partial r  } - 
                       \dfrac{           A_\theta}{         r^2} \right)
          - \sigma u_z b_r
          = \sigma u_z B_r^a
\end{equation}
\begin{equation}\label{eq:2Dt9ge_b}
      \dfrac{1}{\mu}   \dfrac{\partial b_r}{\partial z} - 
                       \dfrac{\partial h_z}{\partial r}
          + \sigma u_z \dfrac{\partial   A_\theta}{\partial z  }
          = \sigma u_z B_r^a
\end{equation}
where $z,~r$ are the axial and radial directions, $B_r^a$ is the radial component
of the applied magnetic field due to the current carrying coil. In this 
problem, the reaction magnetic field has 2-components $b_r$, $b_z$. The 
moving conductor and the air medium has a jump in magnetic permeability. In 
this scenario, the component which is perpendicular to the moving conductor 
$\bf b_\perp$, is continuous across the conductor-air boundary. The parallel
component is discontinuous across the boundary; however, the parallel component
of magnetic field intensity ${\bf h_\parallel}$ is continuous across the 
conductor-air boundary. Thus the variables $b_r$ and $h_z$ are continuous
in this problem. 
By following the 1D case in 
(\ref{eq:wr1Dmc_1}), (\ref{eq:wr1Dmc_2}) the weighted residual
formulation for this 2D problem can be obtained. The (\ref{eq:2Dt9ge}) is
weighted with the Galerkin weight function $N$ and the (\ref{eq:2Dt9ge_b}) 
is weighted with the weight function $\partial N/\partial z$. In addition
to this, there is another variable $h_z$ present in in this formulation. 
For that,
\begin{equation}\label{eq:wr2dhz_b}
    \int_\Omega N \left( \dfrac{\partial A_\theta}{\partial r} +
                         \dfrac{         A_\theta}{         r} \right) -
    \int_\Omega N \mu h_z = 0
\end{equation}
\\
\hspace{3mm}
The third equation (\ref{eq:wr2dhz_b}) is the Galerkin formulation for the 
parallel component $\bf h_\parallel$ and it is expressed as,
\[ \nabla \times {\bf A} \cdot {\bf r_\parallel} - 
             \mu {\bf h} \cdot {\bf r_\parallel} = 0 \]
where, $\bf r_\parallel$ is the tangent vector along the conductor-air
boundary. 
\begin{figure*}
		\centering
		\mbox{\subfloat[]{\label{2dt9br_GL}
		\includegraphics[scale=0.39]{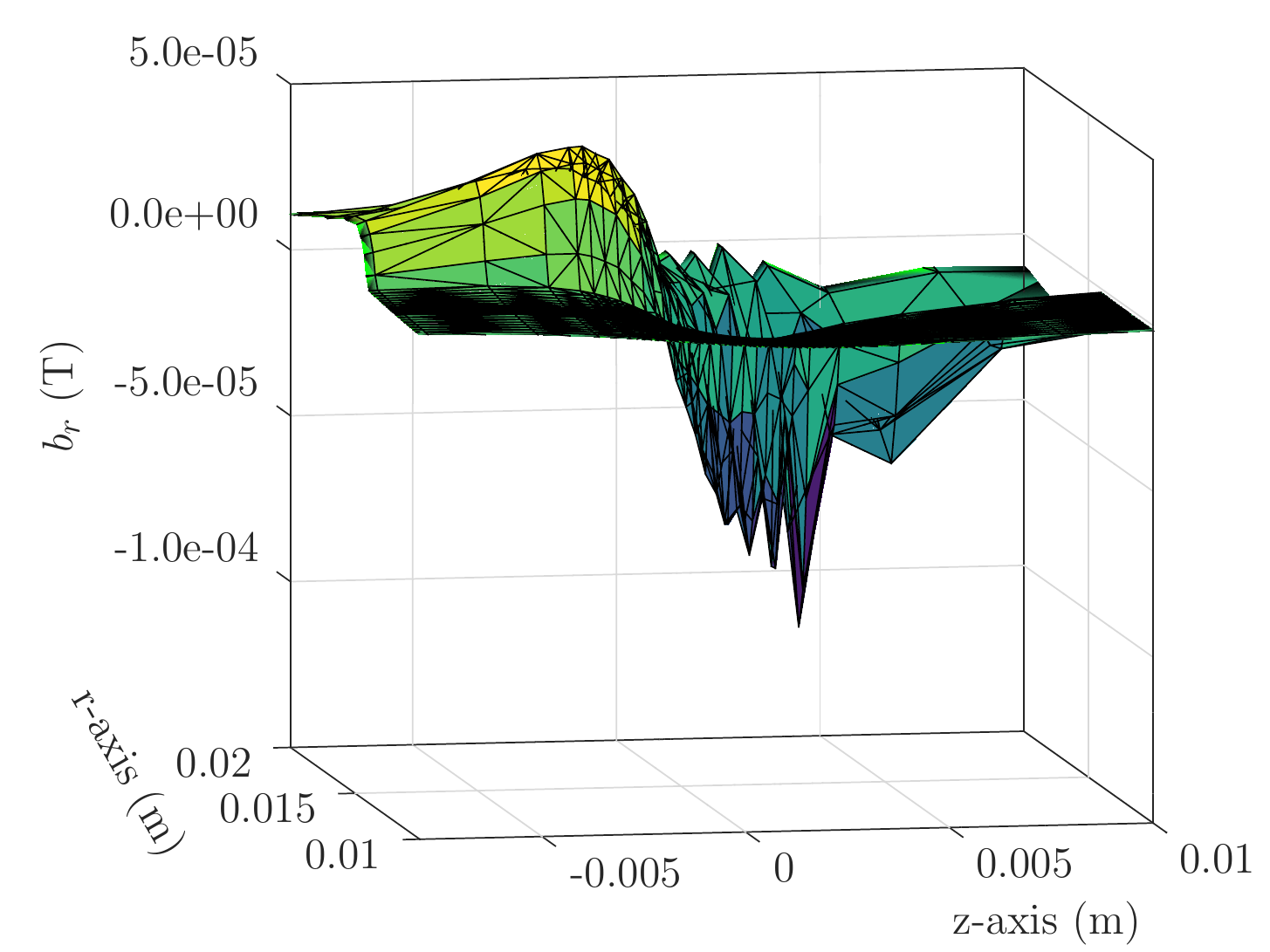}}}
		\mbox{\subfloat[]{\label{2dt9br_SUPG}
		\includegraphics[scale=0.39]{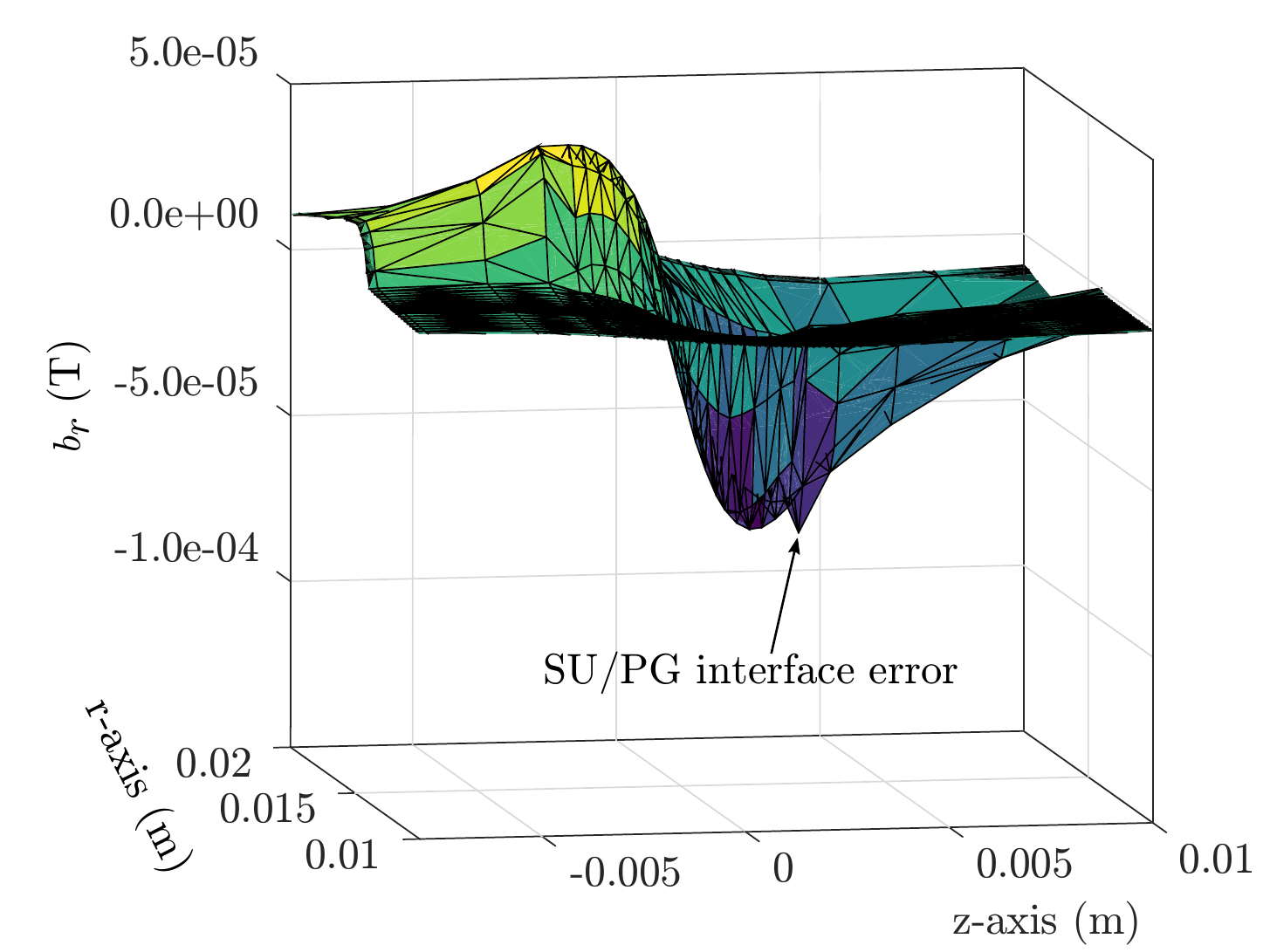}}}
		\mbox{\subfloat[]{\label{2dt9br_WR}
		\includegraphics[scale=0.39]{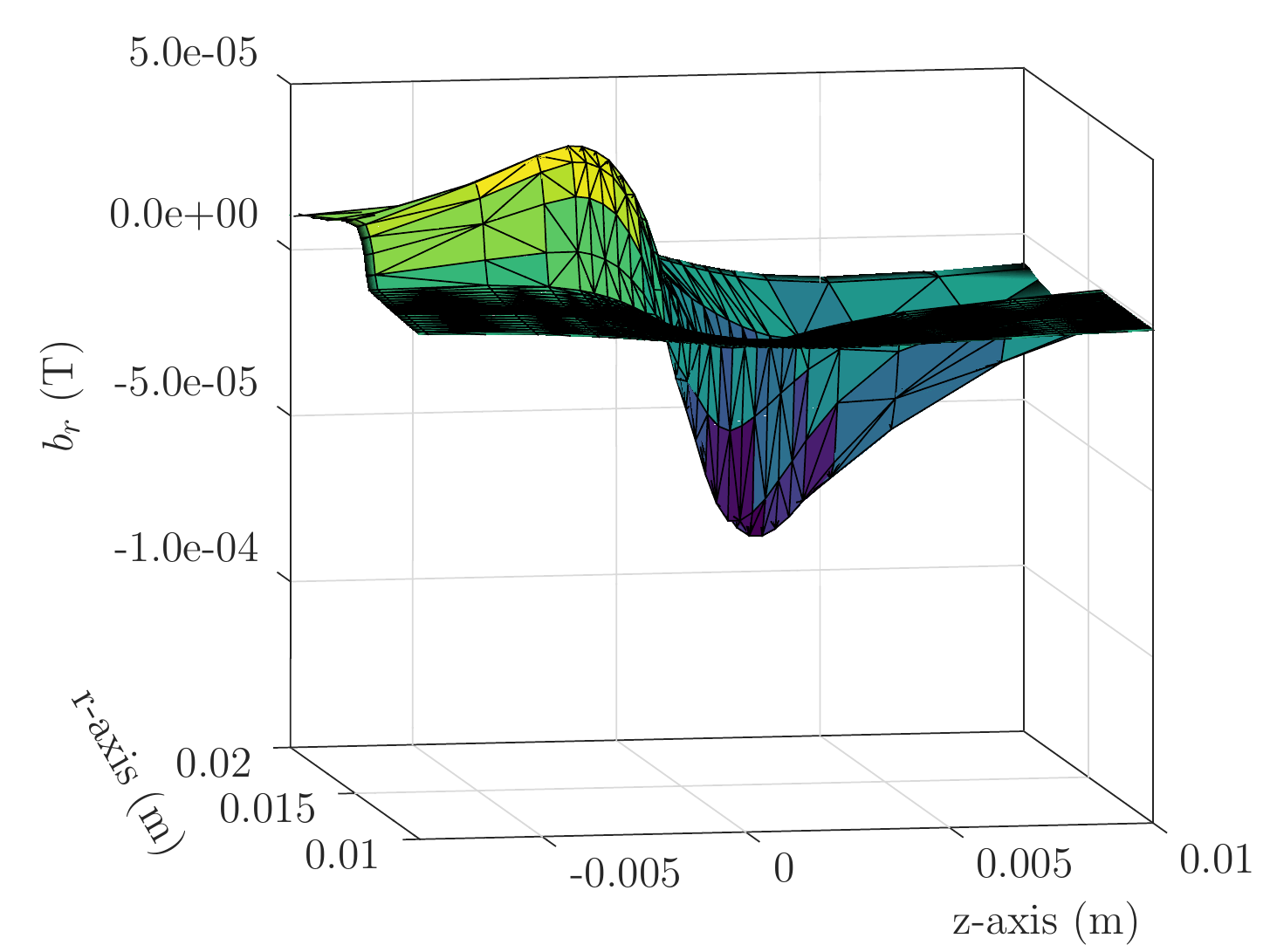}}}
		\caption{Simulation results from the 2D TEAM 9a problem for 
        $u_z$ = 100 $ms^{-1}$ and $\mu_r$ = 50 case
		(a) Galerkin scheme - reaction magnetic field - $b_r$
		(b) SU/PG scheme - reaction magnetic field - $b_r$
        (c) Proposed formulation - reaction magnetic field - $b_r$}
		\label{2dt9br}
\end{figure*}
Simulations are carried out for the TEAM-9a problem, by using the 2280
bilinear quadrilateral elements. The 
results for the $u_z=100ms^{-1}$ and $\mu_r = 50$ case is shown in 
Fig. \ref{2dt9br}. In this, the Fig. \ref{2dt9br_GL} shows the $b_r$ obtained
from the Galerkin scheme and Fig. \ref{2dt9br_WR} shows the $b_r$ obtained
from the weighted residual formulation. It can be seen that the solution from
the proposed formulation is stable. 
{For quick comparison, Fig. {\ref{2dt9br_SUPG}} displays the result obtained
from the SU/PG scheme, showing the error peaking near the interface.}
\begin{figure}
		\centering
		\mbox{\subfloat[]{\label{2dt9_mur1}
		\includegraphics[scale=0.47]{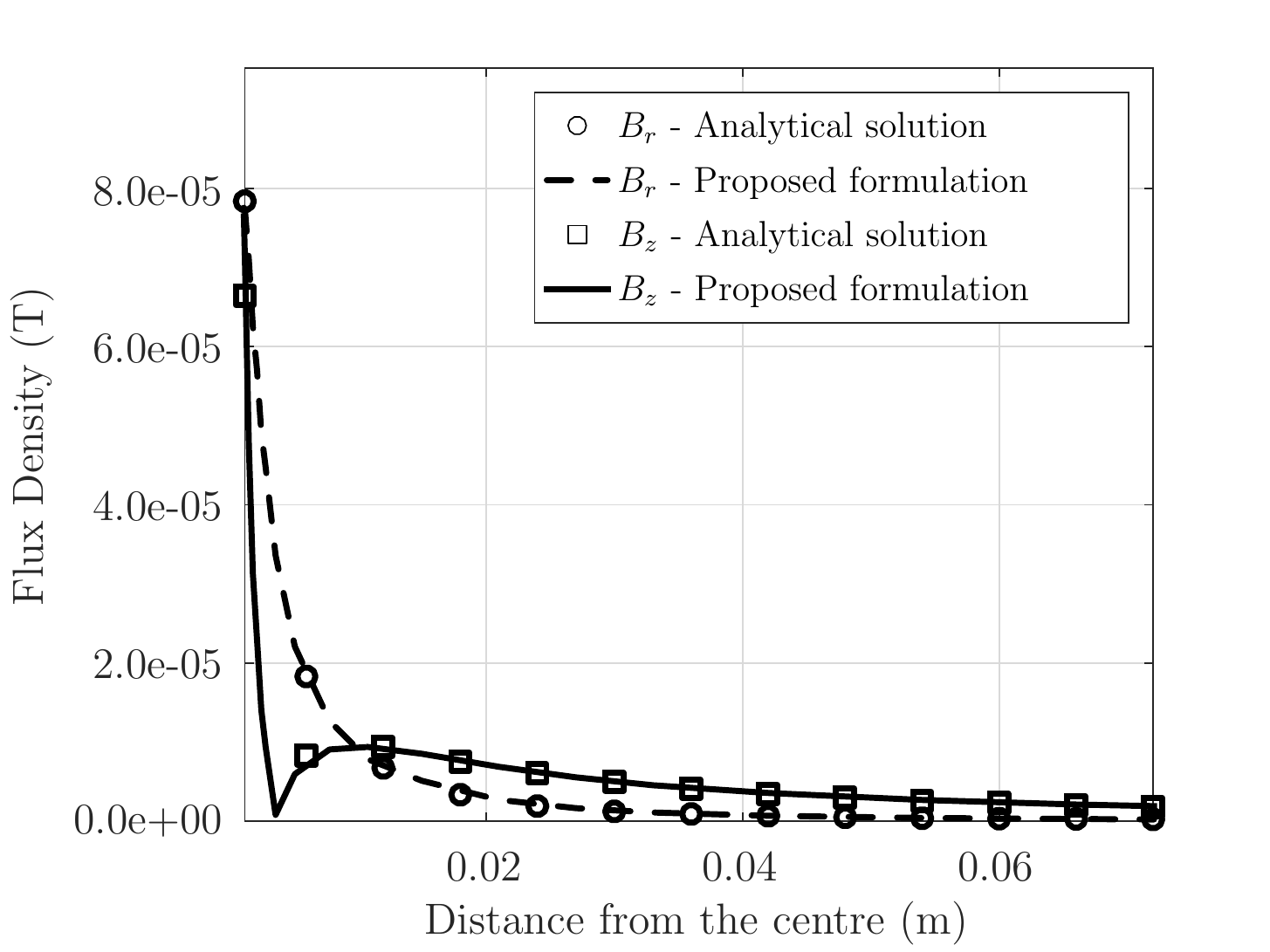}}}
		\mbox{\subfloat[]{\label{2dt9_mur50}
		\includegraphics[scale=0.47]{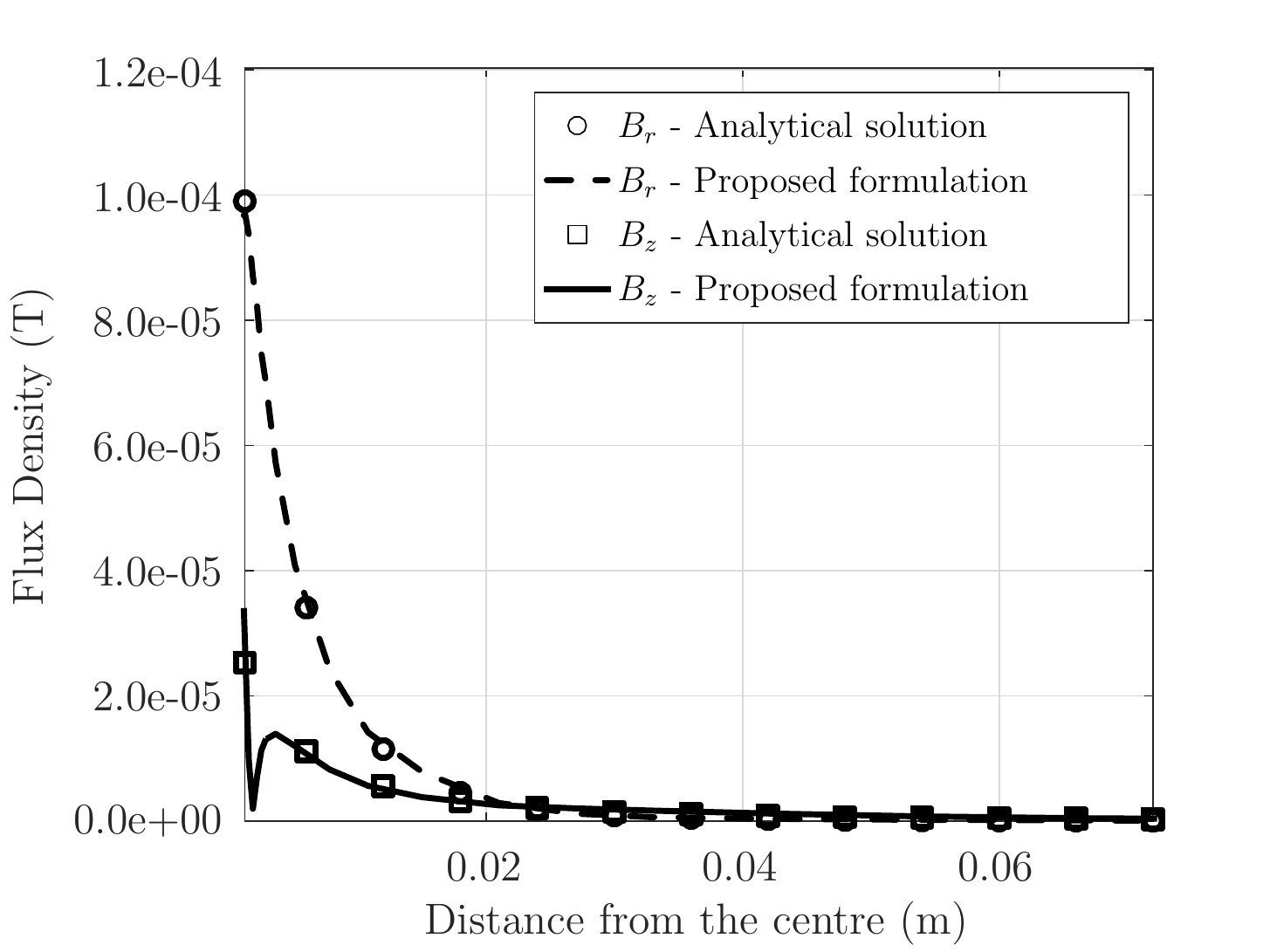}}}
		\caption{Comparison of the total magnetic flux densities from the 
        analytical solution of the 2D TEAM 9a problem and the proposed 
        formulation; for the cases of
        (a) $u_z=100ms^{-1}$, $\mu_r =  1$
        (b) $u_z=100ms^{-1}$, $\mu_r = 50$}
		\label{2dt9_comp}
\end{figure}

The accuracy of the formulation can be seen from Fig. \ref{2dt9_comp},
where the analytical solution of \cite{team9} is compared with the 
solution from the weighted residual formulation. In Fig. \ref{2dt9_mur1}, 
the comparison is made for the non-magnetic conductor 
and the velocity of $u_z = 100 ms^{-1}$; and in Fig. \ref{2dt9_mur50}, 
the comparison is made for the conductor with $\mu_r=50$ and the velocity 
of $u_z = 100 ms^{-1}$.
In the next subsection, the weighted residual formulation for the 
3D case is described.


\section{3D simulation of TEAM-9a problem} \label{sec:4}


The coupled form of the governing equations for the 3D moving conductor problem 
can be written as \cite{sus3, su2},

\begin{equation} \label{eq:3dge2}
\nabla \cdot (\sigma \nabla \phi) - \nabla \cdot (\sigma ~ {\bf{u}} \times \nabla \times {\bf{A}}) =  \nabla \cdot (\sigma ~ \bf{u} \times \bf{B^{a}})
\end{equation}
\begin{equation} \label{eq:3dge1}
\sigma \nabla \phi ~-~  (\nabla \cdot \dfrac{1}{\mu} \nabla) {\bf{A}} - \sigma~ {\bf{u}} \times {\bf{b}} = \sigma~ {\bf{u}} \times {\bf{B^{a}}}
\end{equation}
\begin{equation} \label{eq:3dge1a}
\sigma \nabla \phi ~+~  \nabla \times \dfrac{\bf b}{\mu} - \sigma~ {\bf{u}} \times \nabla \times {\bf{A}} = \sigma~ {\bf{u}} \times {\bf{B^{a}}}
\end{equation}
For the 3D case, the weighted residual formulation is similar to the 2D 
formulations presented above. 
In this, the 
(\ref{eq:3dge2}), (\ref{eq:3dge1}) have the Galerkin scheme with weight
function $N$ and the (\ref{eq:3dge1a}) is weighted so as to have 
\emph{bilinear form} in the formulation. Hence, the perpendicular components, 
${\bf b}_\perp \rightarrow b_x, b_y$ of (\ref{eq:3dge1a}) are weighted with
the weight function $\partial N/\partial z$; and the parallel component,
${\bf h}_\parallel \rightarrow h_z$ is formulated as,
\begin{equation}\label{eq:wr3d_hz}
    \int_\Omega N \left( \dfrac{\partial A_y}{\partial x} -
                         \dfrac{\partial A_x}{\partial y} \right) -
    \int_\Omega N \mu h_z = 0
\end{equation}
which is similar to (\ref{eq:wr2dhz_b}).
\begin{figure}
		\centering
		\mbox{\subfloat[]{\label{team9a_sch}
		\includegraphics[scale=0.32]{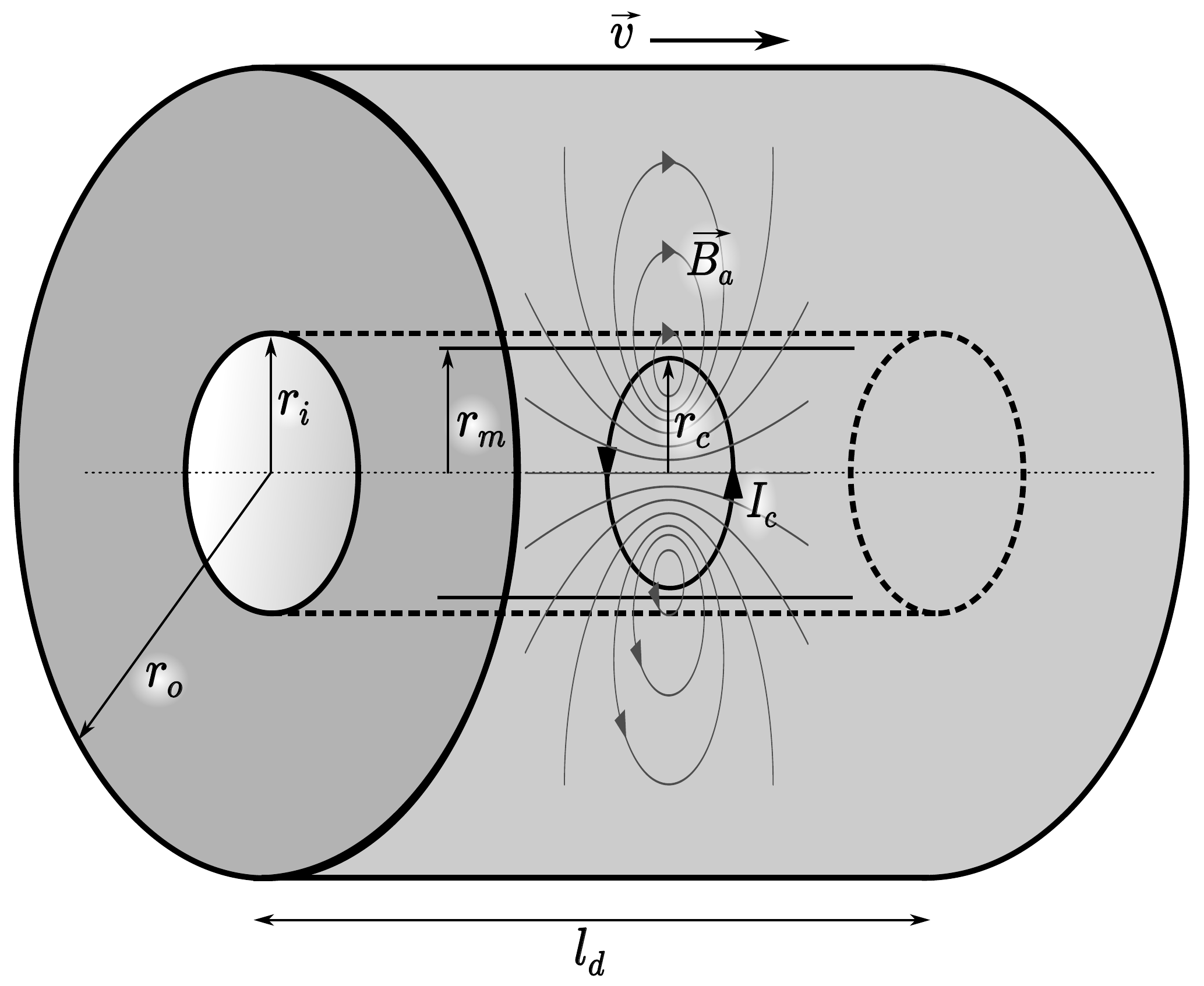}}}
		\mbox{\subfloat[]{\label{team9a_3dmesh}
		\includegraphics[scale=0.18]{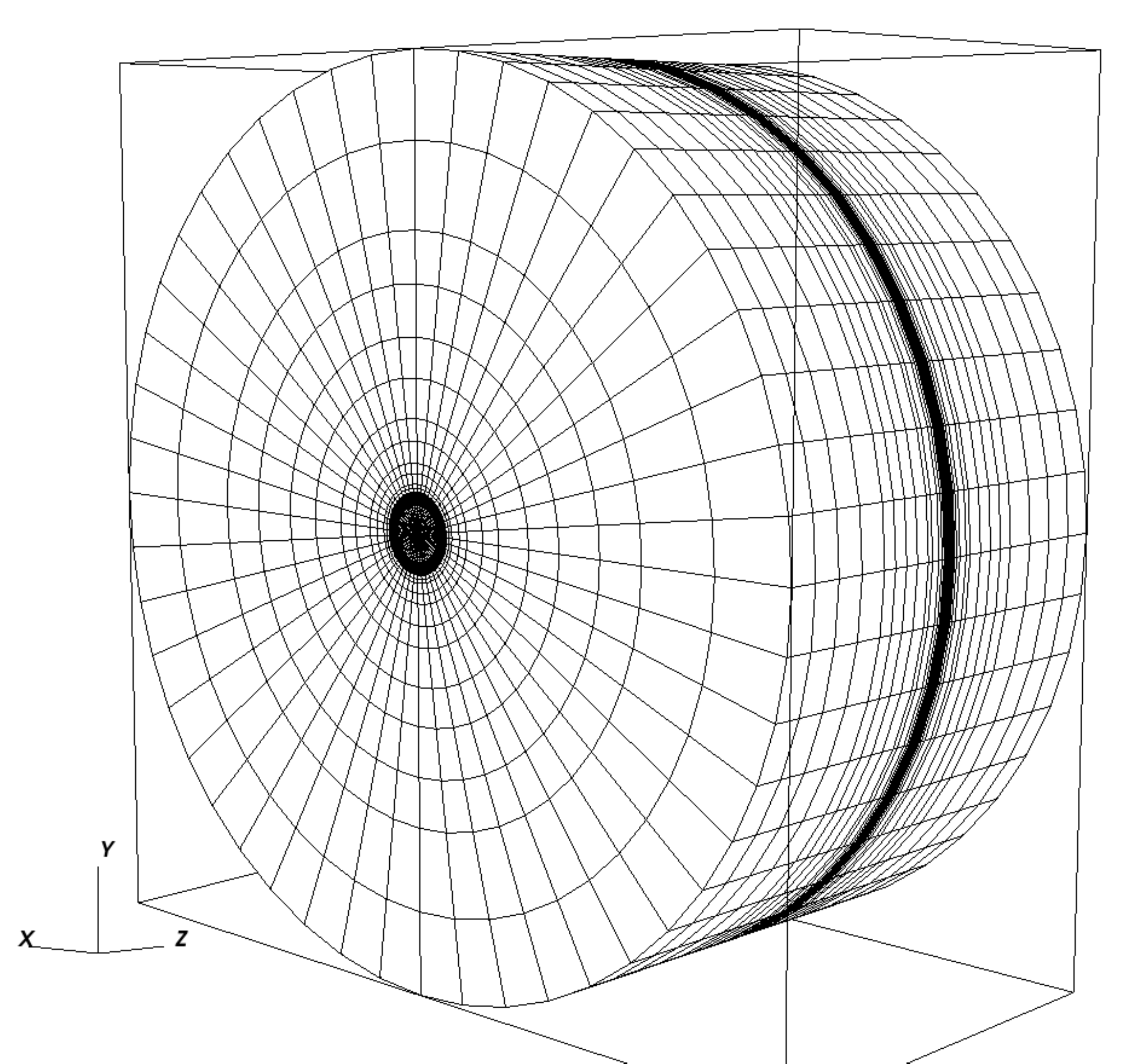}}}
		\caption{Description of the TEAM 9a problem
		(a) Schematic representation of the TEAM 9a moving conductor
		problem.
        (b) Finite element mesh employed}
		\label{3dteam9}
\end{figure}
In order to test the formulation in 3D, the TEAM-9a problem in $x,y,z$-coordinate
system is chosen. In the cartesian coordinates the TEAM-9a problem
loses its symmetry and becomes a 3D test case with the materials having
different permeability and conductivity. The results from the 2D simulation,
can serve as a reference to test the correctness of the solution obtained
from the 3D case. The schematic of the 
problem in 3D is shown in Fig. \ref{team9a_sch}. The simulation is carried out
with 72000 trilinear hexahedron elements and its finite element mesh 
is shown in Fig. \ref{team9a_3dmesh}.
\begin{figure}
		\centering
		\mbox{\subfloat[]{\label{3dt9br_GL}
		\includegraphics[scale=0.50]{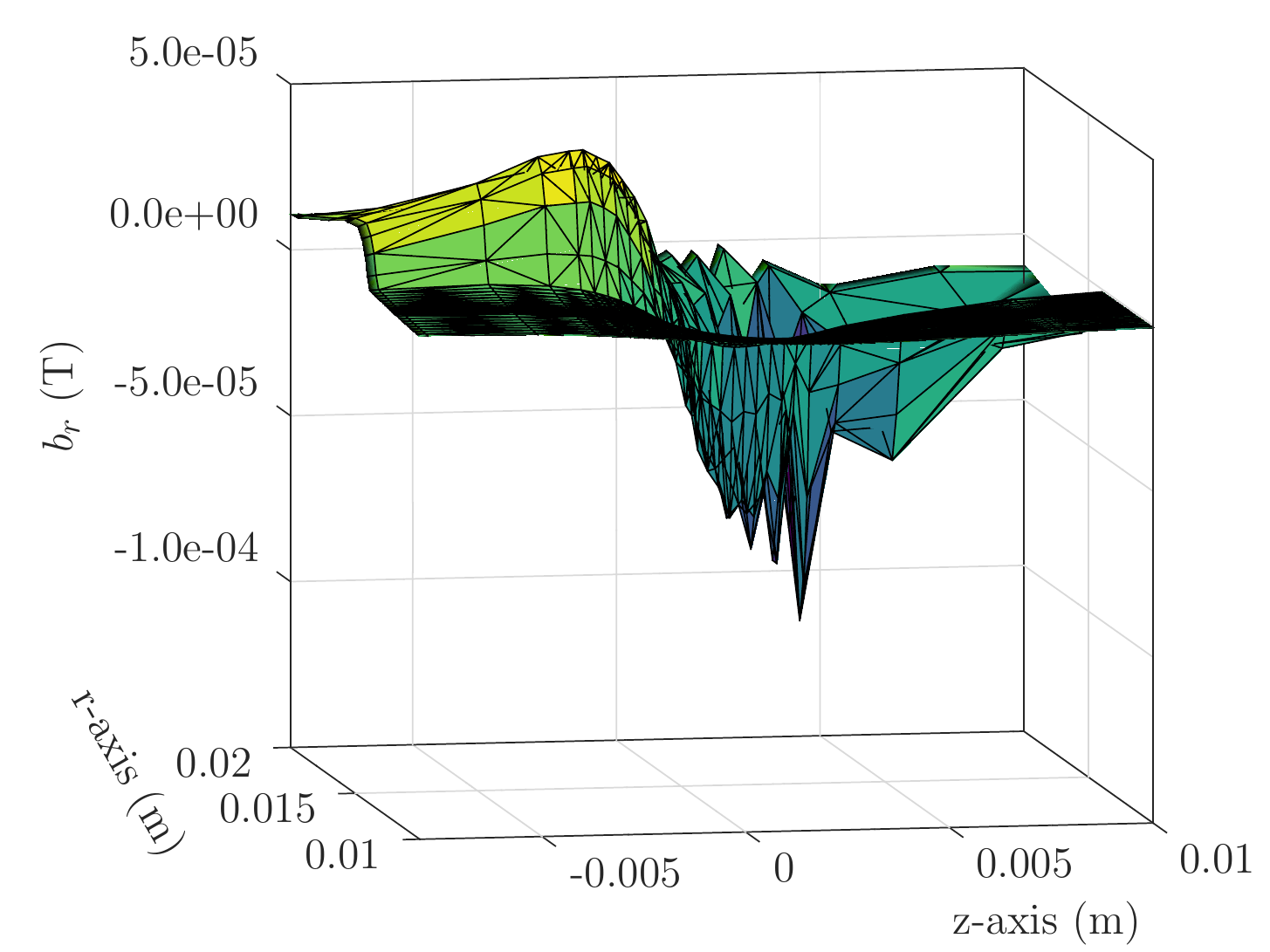}}}
		\mbox{\subfloat[]{\label{3dt9br_WR}
		\includegraphics[scale=0.50]{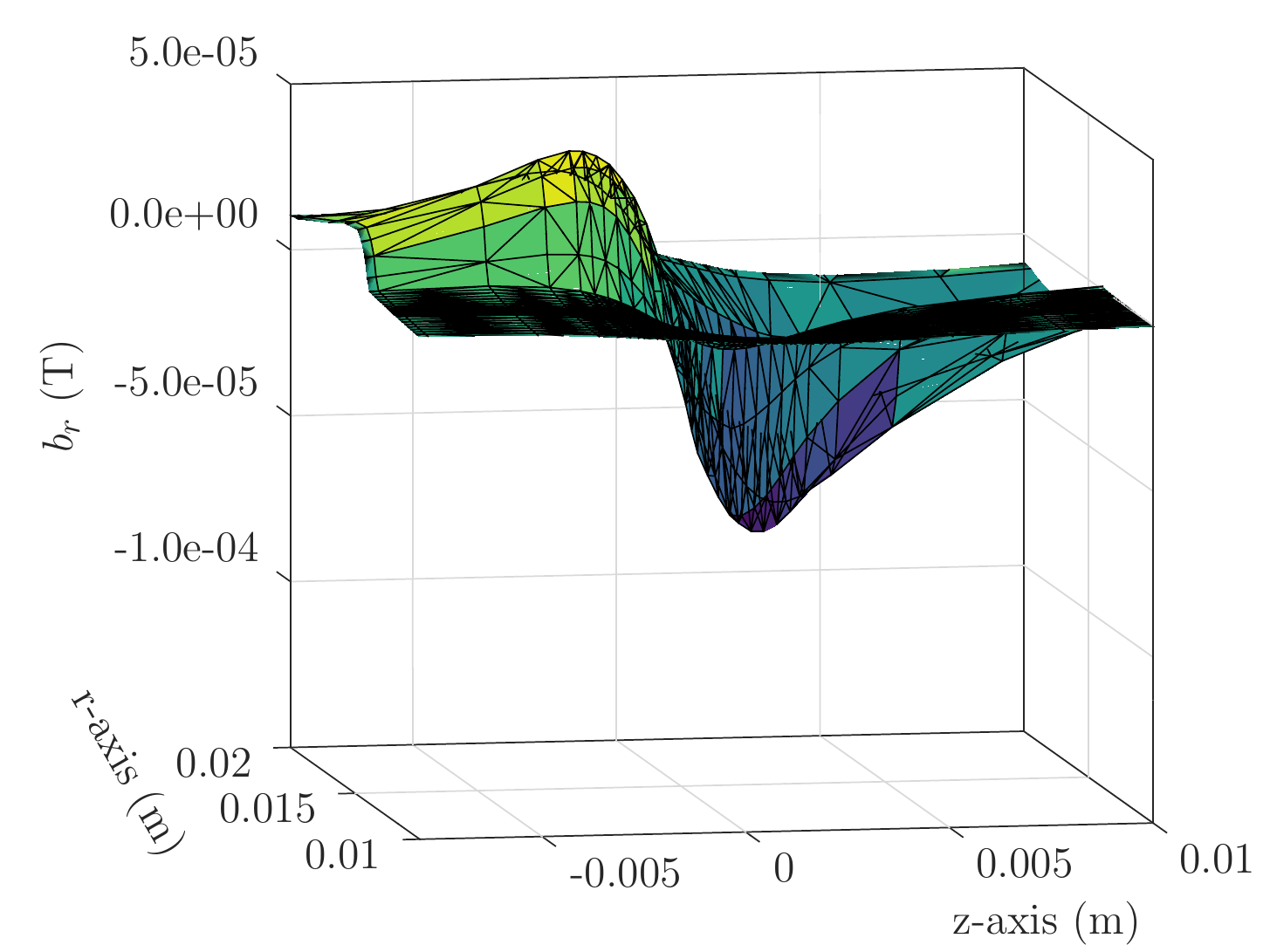}}}
        \caption{Simulation results ($b_r$) from the 3D TEAM 9a problem, at the 
        cross-section along the $\theta \approx 0$ plane, for $u_z$ = 100 
        $ms^{-1}$ and $\mu_r$ = 50 case
		(a) Galerkin scheme
        (b) Proposed formulation}
		\label{3dt9br}
\end{figure}
The reaction magnetic field along the radial direction is plotted in 
Fig. \ref{3dt9br}. The values are taken along the $\theta \approx 0^o$ plane.
When comparing the Fig. \ref{3dt9br} with the Fig. \ref{2dt9br}, it can be
seen that the 3D results exactly resemble the results obtained from the 2D 
case. It may be noted that for the 3D problem, the axial $z$ and radial $r$
discretisations are kept nearly identical to the 2D case. Upon varying the
angular discretisation, it is observed that the results from the 3D simulations
are becoming more accurate with the increasing angular discretisation.
Apart from accuracy, it can also be easily noted that the result from
weighted residual formulation is stable (see Fig. \ref{3dt9br_WR}).
Thus, the proposed formulation performs consistently in 3D as well.


\section{Discussion} \label{sec:5}


The formulation presented here, holds few desired characteristics which are 
absent in the existing schemes. 


i) The formulation is parameter-free. In upwinding schemes, the stability 
and the accuracy of the formulation relies on the correct value of 
stabilization parameter \cite{quada1, fic2}. 
{For the 1D problems, an accurate expression for the stabilization-parameter is
derived by matching the analytical solution. Hence, for the 1D cases,
the SU/PG scheme would perform better than the proposed formulation.
However, for the 2D and 3D cases, the SU/PG scheme is known to
suffer from the error at the transverse boundary} 
\cite{discop, soldreview1, soldreview2, nemosu, sus3}.
To resolve this, iterative techniques are suggested in the literature;
which require a repeated calculation of the FEM solution to arrive at the 
correct solution \cite{fic2, soldreview1, soldreview2}. 
These iterative techniques make the problem non-linear and increases
the computation burden by several times.
{The typical number of iterations is found to be around in $\sim100s$ or more;
some cases also found to be non-converging} \cite{soldreview2}.
The presented formulation does not require a stabilization parameter and hence,
it is free from such a computational burden to arrive at a stable and accurate
solution for the moving conductor problems.

ii) The formulation does not require any special representation of the 
input magnetic field $\bf B^a$, unlike \cite{su1, su2}.
This is because, the numerical stability is not brought in by the 
$pole-zero$ cancellation of the input magnetic field. 
The formulation is inherently stable and the source term can be of any form. 

iii) {It may be noted that the reaction magnetic field, which is the practical 
output of the simulation, is measured from $\nabla \times {\bf A}$ at the 
interior point(s) inside the element; instead of obtaining ${\bf b}$ or ${\bf h}$
from the auxiliary equation. This provides one consistent and
simple way to obtain the magnetic field for cases involving multiple 
magnetic materials. Also, the equations are observed to be numerically coupled,
providing the order of convergence of 1 for both $\nabla \times {\bf A}$
and, ${\bf b}$ or ${\bf h}$.}


\section{Summary and Conclusion} \label{sec:5}

The classical Galerkin finite element method, when applied to moving
conductor problems, is known to lose numerical stability. This is due to
the inability of the central weighted schemes to handle the dominant
first derivative in the governing equation. The common strategy is to 
upwind the formulation, to have more weight along the flow direction.
The correct amount of the upwind is decided by the stabilization
parameter $\tau$. Upwind schemes are known to have other issues and 
various solutions are suggested in the finite element literature, which 
include the iterative solution strategies \cite{discop1, discop, fic2}.

In this work, a different route is taken for the simulation of linear
moving conductor problems. The central weighting is retained
and the first derivative is excluded, by having auxiliary equation(s).
In this way, the formulation is not upwinded and remains parameter-free.
The stability of the formulation is shown in 1D with the help of 
$Z$-transform, as well as, with the numerical examples. In addition to this,
the accuracy and the stability of the formulation is shown with 
the help of different 1D and 2D cases; including the cases with materials
having different conductivity and magnetic permeability. Then the formulation
is verified with a 3D moving conductor simulation; stable and consistent
solutions are observed.








\bibliographystyle{IEEEtran} 
\bibliography{./References/References}


\begin{IEEEbiography}[{\includegraphics[width=1in,height=1.25in,clip,keepaspectratio]{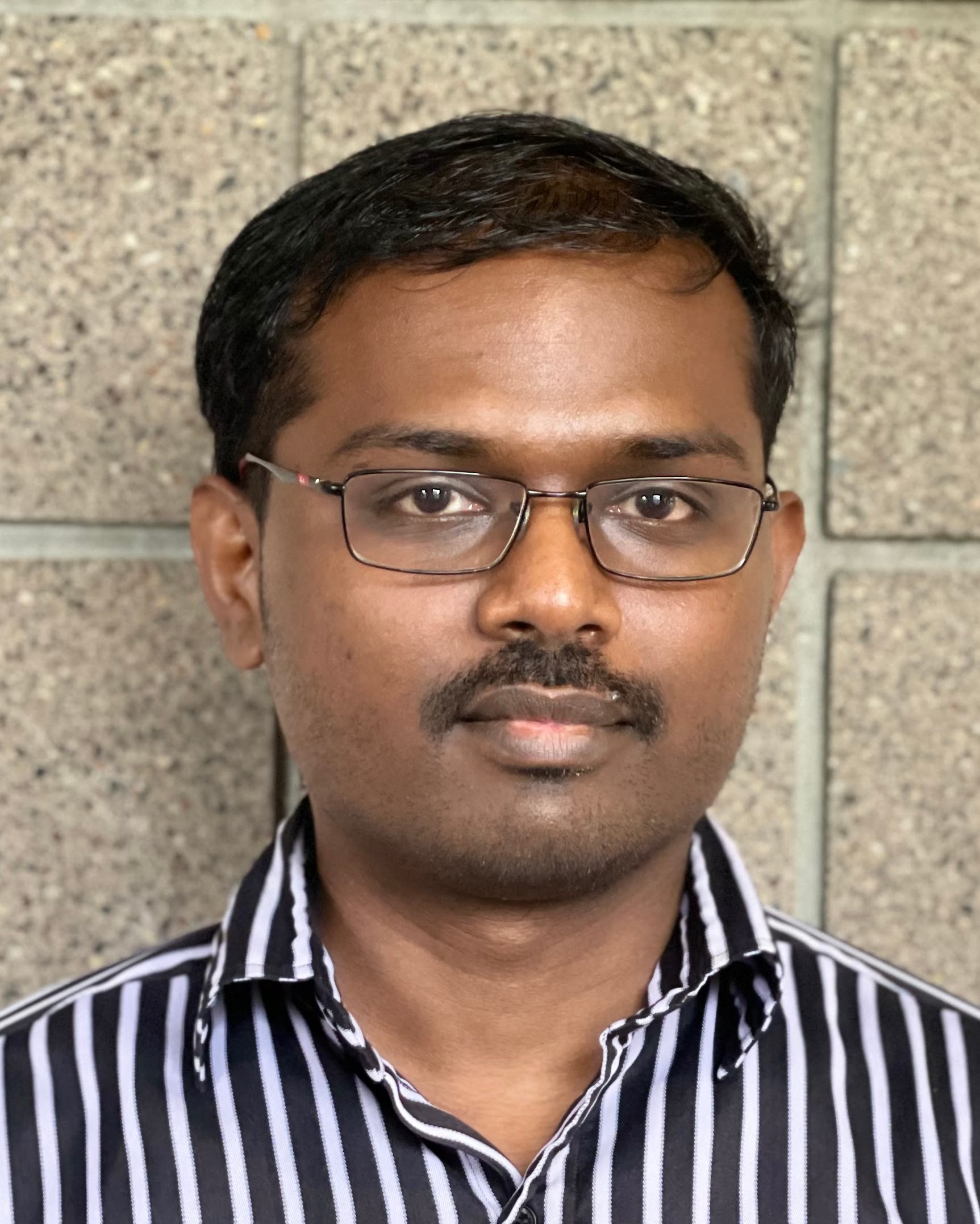}}]{Sethupathy Subramanian}
received the masters and doctrate degrees in electrical engineering from IISc, 
Bangalore, India in 2011 and 2017 respectively. He is currently pursuing his 
graduate research at the Department of Physics and Astronomy, University of 
Notre Dame, USA.

His research interests, pertinent to electrical engineering, 
include computational electromagnetics, numerical stability, 
finite element and edge element methods. 
\end{IEEEbiography}

\begin{IEEEbiography}[{\includegraphics[width=1in,height=1.25in,clip,keepaspectratio]{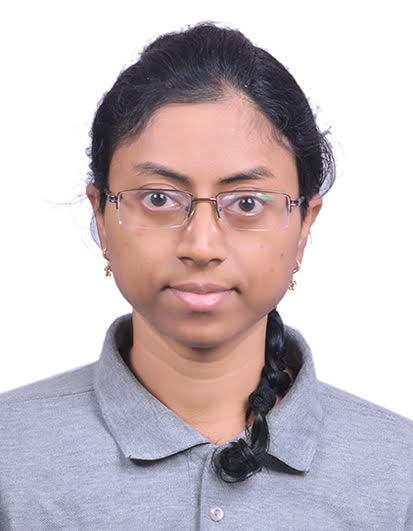}}]{Sujata Bhowmick}
received the B.E. degree in electrical engineering from IIEST, Shibpur, India, in 2006, and the M.E. degree in electrical engineering from the IISc, Bengaluru, India, in 2011. She received the Ph.D. degree from the Department of Electronic Systems Engineering,
IISc, Bengaluru, India, in 2019.

Her research interests include power electronics for renewable resources, single-phase grid-connected power converters, 
computational electromagnetics, finite element and 
edge element methods.
\end{IEEEbiography}
\end{document}